\documentclass[11pt]{article}
\usepackage{amsmath}
\usepackage{amssymb}
\usepackage{amsfonts}
\usepackage{eucal}
\newtheorem{theorem}{Theorem}[section]
\newtheorem{proposition}[theorem]{Proposition}

\newtheorem{corollary}[theorem]{Corollary}

\newtheorem{lemma}[theorem]{Lemma}
\newtheorem{remark}[theorem]{Remark}

\begin{document}
\title{Almost Sure Local Limit Theorem for the Dickman distribution
}\author{Rita Giuliano\thanks{Address: Dipartimento di
Matematica "L. Tonelli", Universit\`a di Pisa, Largo Bruno
Pontecorvo 5, I-56127 Pisa, Italy. e-mail:
\texttt{giuliano@dm.unipi.it}} \and Zbigniew   Szewczak
\thanks{Address: Nicolaus Copernicus University,
Faculty of Mathematics and Computer Science, ul. Chopina 12/18
87-100 Toru\'n, Poland. e-mail: \texttt{zssz@mat.uni.torun.pl}}
\and Michel Weber\thanks{Address: IRMA, UMR 7501, 7, rue Ren\'e--Descartes, 67084 Strasbourg Cedex, France.  e-mail: \texttt{ michel.weber@unistra.fr}}}
\maketitle
\maketitle

\begin{abstract}

\noindent In this paper we present a new correlation inequality and use it for proving an Almost Sure Local Limit Theorem for the so--called Dickman distribution. Several related results are also proved.


\end{abstract}

\section{Introduction}
 The Dickman function plays an important role
in analytic number theory; see \cite{HT1},\cite{T} for further information; see also the very recent paper \cite{DDD}, where a new application is given. 

\medskip
\noindent
Besides its importance in number theory, the Dickman
function also appears in a large number of problems in several other fields: probability, informatics, algebra; we refer to the paper \cite{HT} (where the new example of Hoare's quickselect algorithm is illustrated) and the references therein.

\medskip
\noindent
In the same paper \cite{HT}, a Local Limit Theorem concerning the Dickman function is stated without proof (the authors only refer to Corollary 2.8 in \cite{ABT}; but the use of this result is not at all easy). 

\medskip
\noindent
Whenever a Local Limit Theorem exists, one can wonder whether it can be accompanied by an Almost Sure version: see for instance the papers \cite{GAS}, \cite{GAW} and \cite{W} for examples concerning the Local Limit Theorem and the Almost Sure Local Limit Theorem for partial sums of i.i.d. random variables; see also \cite{GAS1}  for examples of an Almost Sure Local Limit Theorem for Markov chains.

\medskip
\noindent
The purpose of the present paper is twofold: first, we give a detailed proof of the Local Limit Theorem announced in \cite{HT}; second, we answer affirmatively the natural question whether some sort of Almost Sure Local Limit Theorem can be stated and proved. 

\medskip
\noindent
It is worth noting that both our proof of the Local Theorem and the proof of the Almost Sure Local Theorem rely on a new result  (Proposition \ref{stimabase}) that links  the behaviour of the distribution function of the involved partial sums with their local behaviour.

\medskip
\noindent
As it often happens in the theory of Almost Sure Theorems, the crucial point for the proof of our Almost Sure Local Theorem is a new correlation inequality that can have some interest on its own.

\medskip
\noindent
The paper is organized as follows: in section 2 we introduce the notations and present the main results of the paper, i.e. the Local Limit Theorem (Theorem \ref{local}) and the Almost  Sure Local Limit Theorem (Theorem \ref{ASLLT}). In section 3 Proposition \ref{stimabase} and the Local Limit Theorem are proved; we also discuss two important estimations (Propositions \ref{W1} and \ref{W2}) that will be used in the proof of the correlation inequality of section 5; section 4 contains a convergence result (Proposition \ref{ZS}) for the characteristic functions of the partial sums under study, used in the proof of Proposition \ref{stimabase};  section 5 contains the basic correlation inequality; in section 6 we prove a new general form of the Almost Sure Theorem, suitable for us since it takes into account the particular form of the correlation inequality; this result is also interesting in that, loosely speaking, it can be considered as a generalization of a previous one by T. M
 ori (see Theorem 1 in \cite{M}). Section 7 contains the proof of the Almost Sure Local Limit Theorem.  Finally, section 8 is a sort of appendix in which we have stated and proved a new form for the cumulants of the Bernoulli distribution.

\bigskip
\noindent
{\bf Notation}. By the symbol $C$ we denote a positive constant,  the value of which may change from one case to another. We shall not make any distinction between absolute constants and constants depending on some parameter of the problem.

\section{The main results}
Let $\rho$ be the Dickman function, i.e. the function defined on $[0, + \infty )$ by the two conditions
$${\rm (i)}\,\,\rho(x) =1, \quad 0 \leqslant x \leqslant 1;  \qquad\qquad {\rm (ii)}\,\,x \rho^\prime (x) + \rho(x-1)=0, \quad x >1.$$ It is known (see \cite{HT1}, Lemma 2.6) that 
$$\int_0^\infty \rho(x)\, dx = e^\gamma,$$
where $\gamma$ is the Euler--Mascheroni constant; hence $x \mapsto e^{-\gamma}\rho (x) $, $x \geqslant 0$,
is a probability density, known as the {\it Dickman density}. The distribution function with this density is called the {\it Dickman distribution} and will be denoted with $D$. Thus its probability density is $D^\prime (x)=  e^{-\gamma}\rho (x) $. Some properties of the Dickman function and the Dickman distributions will be discussed in the next section.

\bigskip
\noindent
We are interested in the probabilistic model introduced in \cite{HT}: precisely, let $(Z_k)_{k \geqslant 1}$ be independent and such that, for each $k$, 

$$Z_k = \begin{cases}
1 & \hbox {with probability } \frac{1}{k}\\ 0 & \hbox {with probability }  1-\frac{1}{k}.
\end{cases}$$For every pair of integers $(m,n)$ with $0\leqslant m <n$ denote $$T_m^n = \sum_{k=m+1}^n k Z_k.$$
  For simplicity we also put  $T_n=T_0^n $.
  
  \bigskip
  \noindent
  Here are the main results of this paper.
  
\begin{theorem} {\bf(Local Limit Theorem).} \label{local} 
 Let $(\kappa_n)_{n\geqslant 1}$ be any sequence of integers with $\lim_{n \to \infty}\frac{\kappa_n}{n}=x >0$. Then 
 $$\lim_{n \to \infty}nP(T_n = \kappa_n)= e^{- \gamma}\rho(x).$$
 \end{theorem}
 
 \noindent
 The proof of Theorem \ref{local} is in Section 3.
  \begin{theorem} \label{ASLLT}{\bf (Almost Sure Local Limit Theorem).} Let $x>0$ be fixed and let $\boldsymbol \kappa= (\kappa_n)_{n \geqslant 1}$ be a strictly increasing sequence with $\lim_{n \to \infty}\frac{\kappa_n}{n}=x >0$. Then
 $$\lim_{N \to \infty} \frac{1}{\log N}\sum_{n=1}^N1_{\{T_n = \kappa_n\}}= e^{- \gamma}\rho(x), \qquad a.s.$$
 
 \end{theorem}

 \noindent
 The proof of Theorem \ref{ASLLT} is in Section 7.
 
 \bigskip\noindent

 \begin{corollary} We have
 $$\lim_{N \to \infty} \frac{1}{\log N}\sum_{n=1}^N1_{\{T_n = n\}}= e^{- \gamma}, \qquad a.s.$$
As a consequence, for every $x\geqslant 1$,
 $$\lim_{N \to \infty} \frac{\sum_{n=1}^N1_{\{T_n = [xn]\}}}{\sum_{n=1}^N1_{\{T_n = n\}}}= \rho(x), \qquad a.s.$$
 \end{corollary}

\bigskip
 \begin{remark} \rm
 This corollary suggests two simulation procedures
for estimating (i) Euler's constant $\gamma$ and (ii) the values of  Dickman's function $\rho$.


\noindent
Anyway, we have not investigated the goodness of these methods, nor compared them with the existing ones. For the values of the Dickman function see for instance \cite{HT1}, Corollary 2.3. \end{remark}
\section{Preliminaries for the Dickman distribution I: some known facts and a new result}
\noindent 
Let $\rho$ and $D$ be the Dickman function and the Dickman distribution respectively. It is easy to see that
\begin{equation} \label{equazione}
D(x) - D(x-1)=  e^{-\gamma} x\rho (x)=x D^\prime (x), \quad x \geqslant 1.
\end{equation}
In fact, denoting provisorily $f(x)=D(x) - D(x-1)$ and  $ e^{-\gamma} x\rho (x)= g(x)$, $x \geqslant 1$, we have, by (i) of section 2,
$$f(1) = \int_0^1 D^\prime (t) \,dt =  \int_0^1e^{-\gamma}\rho (t)\, dt =e^{-\gamma} =e^{-\gamma} \rho(1) = g(1);$$
and by (ii)
\begin{align*}
&f^\prime(x) =\frac{d}{dx}\Big(\int_{x-1}^x D^\prime (t) \,dt\Big)=\frac{d}{dx}\Big(\int_{x-1}^x e^{-\gamma} \rho (t) \,dt\Big)=e^{-\gamma} \rho(x) -e^{-\gamma} \rho (x-1)\cr & =e^{-\gamma}\Big(\rho(x)+ x \rho^\prime (x)\Big)=\frac{d}{dx}\Big(e^{-\gamma}x\rho(x)\Big)= g^\prime(x), \quad x>1.
\end{align*}
We also recall that the characteristic function of the Dickman distribution is
\begin{equation}\label{characteristicfunction}
\phi(t)=\exp \Big\{\int_0^1 \frac{e^{itu}-1}{u}\Big\}du ,
\end{equation}
see \cite{HT1} again, or \cite{HT}.

\bigskip
\noindent
Let $(Z_k)_{k \geqslant 1}$ and $T_m^n = \sum_{k=m+1}^n k Z_k$ ($0\leqslant m <n$) be as in Section 2.
The characteristic function of $Z_k$ is $$\phi_{Z_k}(t)= 1+ \frac{e^{it}-1}{k}.$$ 
The characteristic function of $T_m^n$ is
 $$\phi_{T_m^n}(t)=\prod_{k=m+1}^n\phi_{Z_k}(tk)=\prod_{k=m+1}^n\Big(1+ \frac{e^{itk}-1}{k}\Big).$$

 \bigskip
 \noindent
 The proof of the following result is identical to the one given in \cite{HT} (Proposition 1) for the case $m_n \equiv 0$.
\begin{proposition}\label{law} Let $(m_n)_{n \geqslant 1}$ be a sequence of integers such that $\lim_{n \to \infty}(n-m_n)= +  \infty$. Then, as $n \to \infty$,
 the sequence $\frac{T_{m_n}^n}{n-m_n}$ converges in distribution to the Dickman law.

 \end{proposition} 
 
 \bigskip
 \noindent Now we present a result that will be crucial for the proof of the correlation inequality of section 5. Its aim is to connect the local behaviour of $T_m^n$ to the behaviour of its distribution function; it can be considered as a quantitative version of Theorem 2.6 in \cite{ABT}.
\begin{proposition}\label{stimabase}
Let $(\kappa_n)_n$ be any increasing sequence of integers. Then, for $n> m \geq 2$, $$\big|(\kappa_n- \kappa_m)P(T_m^n = \kappa_n- \kappa_m)-
 P\big((\kappa_n- \kappa_m)-n < T_m^n \leqslant (\kappa_n- \kappa_m)-(m+1)\big) \big|\leqslant C  \frac{1+\log\frac{n}{m}}{\sqrt{n-m}}.$$
\end{proposition}

\noindent
{\it Proof.} We need a preliminary easy result.
\begin{lemma}\label{iniziale}
Let $T$ be a random variable taking integer values and with characteristic function $\phi_T$. For every integers $\kappa$, $a$ and $b$ with $a <b$ we have the formula
$$P(\kappa- b \leqslant T \leqslant \kappa-a)= \frac{1}{2\pi}\int_{-\pi}^\pi e^{-it \kappa}\Big(\sum_{j =a}^b e^{it j}\Big)\phi_{T}(t)\, dt.$$
\end{lemma}
{\it Proof.} By the inversion formula we can write
\begin{equation*}
P(\kappa- b \leqslant T \leqslant \kappa-a)=\sum_{j=\kappa  -b}^{\kappa -a}P(T=j)=\frac{1}{2\pi}\int_{-\pi}^\pi\Big(\sum_{j=\kappa - b}^{\kappa -a}e^{-it j}\Big)\phi_{T}(t)\, dt,
\end{equation*}
and now, for every $t\in \mathbb{R}$,\begin{equation*}
\sum_{j=\kappa - b}^{\kappa -a}e^{-it j}= e^{-it (k-b)}+ e^{-it (k-b+1)}+ \cdots +e^{-it (k-a)}= e^{-it k}\sum_{j=a}^{b}e^{it j}
\end{equation*}
\hfill$\Box$

We pass to the proof of Proposition \ref{stimabase}. First, by Lemma \ref{iniziale}
\begin{equation}\label{4}
P\big((\kappa_n -\kappa_m)- n < T_m^n \leqslant (\kappa_n -\kappa_m)-(m+1)\big)=\frac{1}{2\pi}\int_{-\pi}^\pi e^{-iu (\kappa_n-\kappa_m)}\Big(\sum_{k =m+1}^n e^{iu k}\Big)\phi_{T_m^n}(u)\, du.
\end{equation}
Moreover, integrating by parts,
 \begin{align}\label{3}
&\nonumber \frac{1}{2\pi i} \int_{-\pi }^{\pi }e^{-iu(\kappa_n -\kappa_m) }\cdot\phi^\prime_{T_m^n}(u)\, du\\ &\nonumber=\frac{1}{2\pi i}\Big\{\phi_{T_m^n}(u)e^{-iu (\kappa_n -\kappa_m) }\Big|_{-\pi } ^{\pi }+ i(\kappa_n -\kappa_m)\int_{-\pi }^{\pi }\phi_{T_m^n}(u)e^{-iu (\kappa_n -\kappa_m)}\, du\Big\}\\&=\nonumber
\frac{1}{2\pi i}\Big\{\Big(\phi_{T_m^n}(\pi )e^{-i\pi (\kappa_n -\kappa_m)}  -\phi_{T_m^n}(-\pi )e^{i\pi (\kappa_n -\kappa_m)}\Big)+i(\kappa_n -\kappa_m)\int_{-\pi }^{\pi }\phi_{T_m^n}(u)e^{-i(\kappa_n -\kappa_m)u }\, du\Big\}\\&=\nonumber\frac{1}{2\pi i}\Big\{\Big(\phi_{T_m^n}(\pi )e^{-i\pi(\kappa_n -\kappa_m) }  -\overline{\phi_{T_m^n}(\pi )e^{-i\pi(\kappa_n -\kappa_m) }}\Big)+2 \pi i (\kappa_n -\kappa_m) P(T_m^n=\kappa_n -\kappa_m)\Big\}\\&=\frac{\mathfrak{Im} \big(\phi_{T_m^n}(\pi )e^{-i\pi(\kappa_n -\kappa_m) } \big)}{\pi}+ (\kappa_n -\kappa_m) P(T_m^n=\kappa_n -\kappa_m)=  (\kappa_n -\kappa_m)P(T_m^n=\kappa_n -\kappa_m),\end{align}
noticing that $\phi_{T_m^n}(\pi )$ and $e^{-i\pi (\kappa_n -\kappa_m)} $ are real (recall that $\kappa_n$ is an integer).
Since
$$\phi^\prime_{T_m^n}(u)=\phi_{T_m^n}(u)\sum_{k =m+1}^n\frac{k\phi_{Z_k}^\prime(ku)}{\phi_{Z_k}(ku)}, $$
subtracting \eqref{4} from \eqref{3} we obtain
\begin{align*}
&(\kappa_n- \kappa_m)P(T_m^n = \kappa_n- \kappa_m)-
 P\big((\kappa_n- \kappa_m)-n < T_m^n \leqslant (\kappa_n- \kappa_m)-(m+1)\big) \\&= \frac{1}{2\pi } \int_{-\pi }^{\pi }e^{-iu(\kappa_n- \kappa_m)}\cdot\phi_{T_m^n}(u)\Big(\sum_{k =m+1}^n\frac{\frac{k}{i}\phi_{Z_k}^\prime(ku)}{\phi_{Z_k}(ku)}  \Big)\, du\\&-\frac{1}{2\pi}\int_{-\pi}^\pi e^{-iu(\kappa_n- \kappa_m) }\Big(\sum_{k =m+1}^n e^{iu k}\Big)\phi_{T_m^n}(u)\, du\\&=
\frac{1}{2\pi } \int_{-\pi }^{\pi }e^{-iu(\kappa_n- \kappa_m) }\cdot\phi_{T_m^n}(u)\Big(\sum_{k =m+1}^n\frac{\frac{k}{i}\phi_{Z_k}^\prime(ku)}{\phi_{Z_k}(ku)}-\sum_{k =m+1}^n e^{iu k}\Big)\, du\\&=\frac{n-m}{2\pi } \int_{-\pi }^{\pi }e^{-iu (\kappa_n- \kappa_m)}\cdot\phi_{T_m^n}(u)\underbrace{\frac{\sum_{k =m+1}^n\frac{\frac{k}{i}\phi_{Z_k}^\prime(ku)}{\phi_{Z_k}(ku)}-\sum_{k =m+1}^n e^{iu k}}{n-m}}_{= \gamma_{m,n} (u)}\, du\\&= \frac{1}{2\pi } \int_{-\pi (n-m) }^{\pi (n-m)}e^{-iu\frac{\kappa_n- \kappa_m}{n-m} }\cdot\phi_{\frac{T_m^n}{n-m}}(u)\,\gamma_{m,n} \Big(\frac{u}{n-m}\Big) \, du.\end{align*}
Hence
\begin{align}\label{intermedia}
&\nonumber\Big|(\kappa_n- \kappa_m)P(T_m^n = \kappa_n- \kappa_m)-
 P\big((\kappa_n- \kappa_m)-n < T_m^n \leqslant (\kappa_n- \kappa_m)-(m+1)\big)\Big|\cr &\leqslant \frac{1}{2\pi }\int_{-\pi (n-m) }^{\pi (n-m)}\Big|e^{-iu\frac{\kappa_n- \kappa_m}{n-m} }\cdot\phi_{\frac{T_m^n}{n-m}}(u)\,\gamma_{m,n} \Big(\frac{u}{n-m}\Big) \Big|\, du \\ &\nonumber\leqslant \frac{1}{2\pi } \Big\{\int_{-\pi (n-m) }^{\pi(n-m) }\Big|e^{-iu\frac{\kappa_n- \kappa_m}{n-m} }\cdot\phi_{\frac{T_m^n}{n-m}}(u)\Big|^2\, du\Big\}^{\frac{1}{2}}\Big\{\int_{-\pi (n-m) }^{\pi (n-m)}\Big|\gamma_{m,n} \Big(\frac{u}{n-m}\Big)|^2\, du\Big\}^{\frac{1}{2}}\\ &\leqslant
\frac{1}{2\pi } \Big\{\int_{-\pi (n-m) }^{\pi (n-m)}\Big|\phi_{\frac{T_m^n}{n-m}}(u)\Big|^2\, du\Big\}^{\frac{1}{2}}\cdot \Big\{2\pi (n-m)\sup_{-\pi \leqslant u\leqslant \pi }|\gamma_{m,n} (u)|^2\Big\}^{\frac{1}{2}}\end{align}
At the end of this proof we shall show that
\begin{equation}\label{bound2}
\sup_{-\pi \leqslant u\leqslant \pi }|\gamma_{m,n} (u)|\leqslant C\frac{1+\log \frac{n}{m}}{n-m} .
\end{equation}
Using \eqref{bound2} in \eqref{intermedia}, we get
\begin{align*}
& \Big\{\int_{-\pi (n-m) }^{\pi (n-m)}\Big|\phi_{\frac{T_m^n}{n-m}}(u)\Big|^2\, du\Big\}^{\frac{1}{2}}\cdot \Big\{2\pi (n-m)\sup_{-\pi \leqslant u\leqslant \pi }|\gamma_{m,n} (u)|^2\Big\}^{\frac{1}{2}}\\&\leqslant C \Big\{\int_{-\pi (n-m) }^{\pi (n-m)}\Big|\phi_{\frac{T_m^n}{n-m}}(u)\Big|^2\, du\Big\}^{\frac{1}{2}}\cdot\frac{1+\log \frac{n}{m}}{\sqrt{n-m}} .
\end{align*}
Since
$$\frac{T_m^n}{n-m}= \frac{T_n}{n}\cdot\frac{n}{n-m}-\frac{\sum_{k=1}^mk Z_k}{n-m} ,$$ putting $W=-\frac{\sum_{k=1}^mk Z_k}{n-m}$ we can write, by independence,
$$\phi_{\frac{T_m^n}{n-m}}(u)= \phi_{\frac{T_n}{n}}\Big(u\frac{n}{n-m}\Big)\cdot \phi_{W}(u)$$
which implies
\begin{align*}
&\Big\{\int_{-\pi (n-m) }^{\pi (n-m)}\Big|\phi_{\frac{T_m^n}{n-m}}(u)\Big|^2\, du\Big\}^{\frac{1}{2}}\cdot\frac{1+\log \frac{n}{m}}{\sqrt{n-m}}\leqslant \Big\{\int_{-\pi (n-m) }^{\pi (n-m)}\Big| \phi_{\frac{T_n}{n}}\Big(u\frac{n}{n-m}\Big)\Big|^2\, du\Big\}^{\frac{1}{2}}\cdot \frac{1+\log \frac{n}{m}}{\sqrt{n-m}}
\\ &=\frac{n-m}{n} \Big\{\int_{-\pi n }^{\pi n}\Big| \phi_{\frac{T_n}{n}}(u)\Big|^2\, du\Big\}^{\frac{1}{2}}\cdot\frac{1+\log \frac{n}{m}}{\sqrt{n-m}}\leqslant\Big\{\int_{-\pi n }^{\pi n}\Big| \phi_{\frac{T_n}{n}}(u)\Big|^2\, du\Big\}^{\frac{1}{2}}\cdot\frac{1+\log \frac{n}{m}}{\sqrt{n-m}}. \end{align*}

\noindent
In Proposition \ref{ZS} of the next section we shall prove that
\begin{equation*}
\int_{-\pi n }^{\pi n}\Big|\phi_{\frac{T_n}{n}}(u)\Big|^2\, du\to \int_{-\infty }^\infty|\phi(u)|^2\,du< \infty,
\end{equation*}
where $\phi$ is as in \eqref{characteristicfunction}; this concludes the proof.

\bigskip\noindent
It remains to prove \eqref{bound2}.
Write
\begin{align*}
&(n-m) \gamma_{m,n}(u) =\sum_{k =m+1}^n\Big(\frac{k}{i}\frac{\phi_{Z_k}^\prime(ku)}{\phi_{Z_k}(ku)}- e^{iu k}\Big)\\&=\sum_{k =m+1}^n\Big(\frac{k}{i}\frac{\frac{ie^{iuk}}{k}}{1+ \frac{e^{iuk}-1}{k}}- e^{iu k}\Big)=\sum_{k =m+1}^n e^{iu k}\Big(\frac{k}{k-1+ e^{iu k}}-1\Big)\\&=\sum_{k =m+1}^n \frac{e^{iu k}(1- e^{iu k})}{k-1+ e^{iu k}}=\sum_{k =m+1}^n \frac{e^{iu k}(1- e^{iu k})(k-1+ e^{-iu k})}{|k-1+ e^{iu k}|^2}.
\end{align*}
Hence

\begin{align*}
(n-m) |\gamma_n(u)| \leqslant 2 \sum_{k =m+1}^n\frac{k}{|k-1+ e^{iu k}|^2}\leqslant 2 \sum_{k =m+1}^n\frac{k}{(k-2)^2}\leqslant C \Big(1+\log \frac{n}{m}\Big),
\end{align*}
since
\begin{align*}
&\big|k-1+ e^{iu k}\big|\geqslant \Big|(k-1)-\big|e^{iuk}\big|\Big|=k-2, \qquad k \geqslant m+1 \geqslant 3\\&
\big|k-1+ e^{-iu k}\big|\leqslant(k-1) +\big|e^{-iuk}\big| =k,\qquad \,\quad  k \geqslant 1
\end{align*}
and
\begin{align*}&\sum_{k =m+1}^n\frac{k}{(k-2)^2}=\sum_{k =m-1}^{n-2}\frac{1}{k}+\sum_{k =m-1}^{n-2}\frac{2}{k^2}\leqslant 
\sum_{k =m-1}^{n-2}\frac{1}{k}+C\leqslant \frac{1}{m-1}+\int_{m-1}^{n-2}\frac{1}{x}\,dx+C\\& \leqslant C +\log \frac{n-2}{m-1}=C +\log \frac{\frac{n}{m}-\frac{2}{m}}{1-\frac{1}{m}}\leqslant C +\log 2\frac {n}{m}=C +\log \frac {n}{m}\leqslant C \Big(1+\log \frac{n}{m}\Big).\end{align*}
\hfill $\Box$

\bigskip
\noindent
Now we can give the

\bigskip
\noindent
{\it Proof of Theorem \ref{local}}.
Assume that we are able to prove that 
$$\lim_{n \to \infty}nP( T_2^n =\widetilde \kappa_n)= e^{- \gamma}\rho(x)$$
for every sequence  of integers $(\widetilde\kappa_n)_{n\geqslant 1}$ with $\lim_{n \to \infty}\frac{\widetilde\kappa_n}{n}=x$. Denote $U= Z_1 + 2Z_2$ and notice that $U$  is independent on ${T}_2^n$ and takes the values 0,1,2,3. Now for each $h =0,1,2,3$, take the sequence $(\widetilde\kappa_n^{(h)})_{n\geqslant 1}$ defined by $\widetilde\kappa_n^{(h)}=\kappa_n-h. $ Since $\frac{\widetilde\kappa_n^{(h)}}{n}=\frac{\kappa_n-h}{n}\to x$ as $n \to \infty$, we have
\begin{align*}
&nP(T_n=\kappa_n )= \sum_{h=0}^3 P(U=h)\big\{nP({T}_2^n= \kappa_n-h)\big\}\\&=\sum_{h=0}^3 P(U=h)\big\{nP({T}_2^n= \widetilde\kappa_n^{(h)})\big\}\to \Big(\sum_{h=0}^3 P(U=h)\Big)e^{- \gamma}\rho(x)=e^{- \gamma}\rho(x).
\end{align*}
and the claim is proved.
So, let $(\widetilde\kappa_n)_n$ be a sequence with $\lim_{n \to \infty}\frac{\widetilde\kappa_n}{n}=x$ . By Proposition \ref{law}
and the continuity of the Dickman distribution we have
\begin{align*}
&\frac{n}{\widetilde\kappa_n}\cdot\big\{P(\widetilde\kappa_n - n < T_2^n \leqslant \widetilde\kappa_n-3)\big\}\\&=\frac{n}{\widetilde\kappa_n}\cdot\big\{P(\frac{\widetilde\kappa_n-n}{n-2}  < \frac  {T_2^n}{n-2}\leqslant \frac{\widetilde\kappa_n- 3}{n-2}   )\big\}\mathop{\to}_n
\frac{1}{x}\big\{D(x) -D(x-1)\big\}= e^{-\gamma}\rho(x),\end{align*}
by \eqref{equazione}. Consider the sequence $(\kappa^{\prime}_n)_n$ defined as $$\kappa^{\prime}_n = \begin{cases}
\widetilde\kappa_n& n \geqslant 3\\ 0 & n=1,2.
\end{cases}$$The estimation of Proposition \ref{stimabase} gives, for $n \geqslant 3$
\begin{align*}
&\big|\widetilde\kappa_nP(T_2^n=\widetilde\kappa_n )-P(\widetilde\kappa_n - n < T_2^n \leqslant\widetilde \kappa_n-3)\big|\\ &=
\big|(\kappa^{\prime}_n -\kappa^{\prime}_2)P(T_2^n=\kappa^{\prime}_n -\kappa^{\prime}_2)-P(\kappa^{\prime}_n -\kappa^{\prime}_2- n < T_2^n \leqslant \kappa^{\prime}_n -\kappa^{\prime}_2 -3)\big|\leqslant C \frac{1+\log \frac{n}{2}}{\sqrt{n-2}},
\end{align*}
and the result follows.
\hfill $\Box$
 \begin{remark}\label{cumulants}
\rm Concerning the proof of Proposition \ref{stimabase}, notice that
\begin{equation}\label{a1}
\frac{\phi_{Z_k}^\prime(t)}{\phi_{Z_k}(t)}= \psi^\prime_{Z_k}(t),
\end{equation}
where $\psi_{Z_k}(t)= \log \phi_{Z_k}(t) $, i.e. the second  characteristic function of $Z_k$.

\bigskip
\noindent
Write
$$\psi_{Z_k}(t) = \sum_{j=1}^\infty c_j^{(k)}\frac{(it)^j}{j!},$$
where $\big(c_j^{(k)}\big)_j$ is the sequence of the cumulants of the $\mathcal{B}\big(1, \frac{1}{k}\big)$ distribution. 
Hence
\begin{equation}\label{b}
\psi^\prime_{Z_k}(t) = \sum_{j=1}^\infty ic_j^{(k)}\frac{(it)^{j-1}}{{(j-1)}!}.
\end{equation}
Denote by $\psi_{\Pi_k}(t)$ the second characteristic function of the Poisson law with parameter $\frac{1}{k}$., i.e.
\begin{equation}\label{c}
\psi_{\Pi_k}(t)= \frac{e^{it}-1}{k}, \qquad \psi^\prime_{\Pi_k}(t)= \frac {i}{k}e^{it}=\frac {i}{k}\sum_{j=1}^\infty \frac{(it)^{j-1}}{(j-1)!}.
\end{equation}
Hence, by \eqref{a1} and \eqref{c},
\begin{align*}
&\gamma_{m,n}(t) = \frac{\sum_{k =m+1}^n\frac{\frac{k}{i}\phi_{Z_k}^\prime(tk)}{\phi_{Z_k}(tk)}-\sum_{k =m+1}^n e^{it k}}{n-m}=
\frac{\sum_{k =m+1}^n\frac{k}{i}\big(\frac{\phi_{Z_k}^\prime(tk)}{\phi_{Z_k}(tk)}- \frac{i}{k}e^{itk }\big)}{n-m}\\&= \frac{\sum_{k =m+1}^n\frac{k}{i}\big(\psi_{Z_k}^\prime(tk)-\psi^\prime_{\Pi_k}(tk) \big)}{n-m}.
\end{align*}
Since, by \eqref{b} and \eqref{c},
\begin{align*}
\psi_{Z_k}^\prime(t)-\psi^\prime_{\Pi_k}(t)= \sum_{j=1}^\infty \Big(ic_j^{(k)}- \frac{i}{k}\Big)\frac{(it)^{j-1}}{(j-1)!},
\end{align*}
we get
\begin{equation}\label{alfa}
\gamma_{m,n}(t) =\frac{\sum_{k =m+1}^n\sum_{j=1}^\infty \Big(k c_j^{(k)}- 1\Big)\frac{(itk)^{j-1}}{(j-1)!}}{n-m}=\sum_{j=1}^\infty\frac{(it)^{j-1}}{(j-1)!} \Big\{\frac{\sum_{k =m+1}^nk^{j-1}(k c_j^{(k)}- 1)}{n-m}\Big\}. \end{equation}
Putting $$\alpha_j^{(m,n)}= {\frac{\sum_{k =m+1}^nk^{j-1}(k c_j^{(k)}- 1)}{n-m}}, $$
we obtain the formula

\begin{equation}\label{gamma}
\gamma_{m,n}(t) =\sum_{j=1}^\infty\frac{(it)^{j-1}}{(j-1)!}\alpha_j^{(m,n)}.
\end{equation}

\bigskip
\noindent
Let $\mathbb{B}(1,p)$ be the Bernoullian law with parameter $p\in (0,1)$ an $c_j(p)$ the $j-$th cumulant of $\mathbb{B}(1,p)$. In Section 7 we give an explicit form for the quantity
$$\frac{c_j(p)}{p}-1;$$
if $p=\frac{1}{k}$, this quantity is precisely the expression $k c_j^{(k)}- 1$ in the previous calculations (see \eqref{alfa}).

\noindent
We believe that the explicit formula of Section 7 can be used for getting  good approximations of $k c_j^{(k)}- 1$, and in turn of $\gamma_{m,n}$ (see \eqref{gamma}).
\end{remark}

\noindent
The following result specifies Proposition 1 of \cite{HT} quantitatively in terms of the characteristic functions.

\begin{proposition}\label{W1}
There exists an absolute constant $C$ such that for all integers $n > m \geqslant 2$ and all real numbers $t$,
$$\bigg|\phi_{\frac{T_m^n}{n-m}}(t)-\exp\Big\{\int_0^1 \frac{e^{itu}-1}{u}\,du\Big\}\bigg|\leqslant f_{m,n}(t),$$ where
$$ f_{m,n}(t)=\exp\Bigg\{Ct^2\Big(\frac{\log\frac{n}{m}}{(n-m)^2}+\frac{m+2}{n-m}\Big)\Bigg\}-1 .$$
\end{proposition}
{\it Proof.} First
\begin{align*}
&\bigg|\exp\Big\{\int_0^1 \frac{e^{itu}-1}{u}\,du\Big\}-\phi_{\frac{T_m^n}{n-m}}(t)\bigg|\\&=\bigg|\phi_{\frac{T_m^n}{n-m}}(t)\bigg|\cdot\Big|\exp\Big\{\int_0^1 \frac{e^{itu}-1}{u}\,du-\sum_{k=m+1}^n\log \Big(1+\frac{e^{it\frac{k}{n-m}}-1}{k} \Big)\Big\}-1\Big|\\&\leqslant \Big|\exp\Big\{\sum_{k=m+1}^n\Big[\int_{\frac{k-1-m}{n-m}}^{\frac{k-m}{n-m}} \frac{e^{itu}-1}{u}\,du-\log \Big(1+\frac{e^{it\frac{k}{n-m}}-1}{k} \Big)\Big]\Big\}-1\Big|\\& \leqslant\exp\Bigg\{\sum_{k=m+1}^n\Big|\int_{\frac{k-1-m}{n-m}}^{\frac{k-m}{n-m}} \frac{e^{itu}-1}{u}\,du-\log \Big(1+\frac{e^{it\frac{k}{n-m}}-1}{k} \Big)\Big|\Bigg\}-1\\&\leqslant \exp\Bigg \{\underbrace{\sum_{k=m+1}^n\Big|\frac{e^{it\frac{k}{n-m}}-1}{k}-\log \Big(1+\frac{e^{it\frac{k}{n-m}}-1}{k} \Big)\Big|}_{\rm (a)}+\\&+\sum_{k=m+1}^n\underbrace{\Big|\int_{\frac{k-1-m}{n-m}}^{\frac{k-m}{n-m}} \frac{e^{itu}-1}{u}\,du-\frac{e^{it\frac{k}{n-m}}-1}{k}\Big|}_{\rm (b)}\Bigg\}-1.
\end{align*}
We shall prove that
$$\sum_{k=m+1}^n\Big|\log\Big(1+\frac{e^{it\frac{k}{n-m}}-1}{k}\Big) -\frac{e^{it\frac{k}{n-m}}-1}{k}\Big|\leqslant Ct^2\Big(\frac{\log\frac{n}{m}}{(n-m)^2}+ \frac{1}{n-m} \Big)\eqno \rm (a)$$
and
$$\Bigg|\int_{\frac{k-1-m}{n-m}}^{\frac{k-m}{n-m}} \frac{e^{itu}-1}{u}\,du-\frac{e^{it\frac{k}{n-m}}-1}{k}\Bigg|\leqslant \frac{(m+1)t^2}{2(n-m)^2}, \qquad\qquad\forall \, k \in (m,n]. \eqno \rm (b)$$
These inequalities give
\begin{align*}
&\bigg|\exp\Big\{\int_0^1 \frac{e^{itu}-1}{u}\Big\}\,du-\phi_{\frac{T_m^n}{n-m}}(t)\bigg|
\leqslant \exp \Bigg\{ Ct^2\Big(\frac{\log\frac{n}{m}}{(n-m)^2}+ \frac{1}{n-m}+\sum_{k=m+1}^n\frac{(m+1)}{2(n-m)^2}\Big)\Bigg\}-1\\&\leqslant\exp \Bigg\{Ct^2\Big(\frac{\log\frac{n}{m}}{(n-m)^2}+\frac{m+2}{n-m}\Big)\Bigg\} -1=f_{m,n}(t).
\end{align*}
(a)
The inequality\begin{equation*}\big|e^{ix}-1\big|= 2 \Big|\sin \frac{x}{2}\Big|\leqslant 2\wedge |x| 
\end{equation*}
applied to $x=\frac{tk}{n-m}$ gives
\begin{equation}\label{firstbound}
\Big|\frac{e^{it\frac{k}{n-m}}-1}{k}\Big|\leqslant \frac{2}{k}\wedge  \frac{|t|}{n-m}.
\end{equation}
From 
\begin{equation*}
|\log (1+u)-u|\leqslant \sum_{j=2}^\infty \frac{|u|^{j}}{j}= |u|^2\sum_{j=0}^\infty \frac{|u|^j}{j+2}, \qquad |u|< 1,
\end{equation*}
applied to $u =\frac{e^{it\frac{k}{n-m}}-1}{k}$ (with $k \geqslant m+1$; notice that $\big|\frac{e^{it\frac{k}{n-m}}-1}{k}\big|<1$ for $k \geqslant m+1\geqslant 3$) we get
\begin{align*}
&\sum_{k=m+1}^n\Big|\log\Big(1+\frac{e^{it\frac{k}{n-m}}-1}{k}\Big) -\frac{e^{it\frac{k}{n-m}}-1}{k}\Big|\\&\leqslant \Big\{\max_{m+1 \leqslant k \leqslant n}\Big|\frac{e^{it\frac{k}{n-m}}-1}{k}\Big|^2\Big\}\sum_{j=0}^\infty \frac{1}{j+2}
\sum_{k=m+1}^n\Big|\frac{e^{it\frac{k}{n-m}}-1}{k}\Big|^j\\& \underbrace{\leqslant}_{by \eqref{firstbound}}  \frac{t^2}{(n-m)^2}\Big(\frac{1}{2}(n-m)+ \frac{2}{3}\sum_{k=m+1}^n \frac{1}{k}+\sum_{j=2}^\infty \frac{1}{j+2}
\sum_{k=m+1}^n\Big(\frac{2}{k}\Big)^j\Big)\\& \leqslant C\cdot  \frac{t^2}{(n-m)^2}\Big((n-m) + \log \frac{n}{m}+ \sum_{j=2}^\infty \frac{2^j}{j+2}\int_2^\infty \frac{1}{x^j}{\rm d}x\Big)\\&\leqslant  C\cdot  \frac{t^2}{(n-m)^2}\Big((n-m) + \log \frac{n}{m}+ 1\Big)  \leqslant Ct^2\Big(\frac{\log\frac{n}{m}}{(n-m)^2}+ \frac{1}{n-m}\Big) . \end{align*}
(b) Putting
 $$\eta_t(x) = \frac{e^{itx}-1}{x}, \qquad x>0, \quad -\pi \leqslant t \leqslant \pi,$$
 we can write
 \begin{align*}
 &\Big|\int_{\frac{k-1-m}{n-m}}^{\frac{k-m}{n-m}} \frac{e^{itu}-1}{u}\,du-\frac{e^{it\frac{k}{n-m}}-1}{k}\Big|=
 \Big|\int_{\frac{k-1-m}{n-m}}^{\frac{k-m}{n-m}} \frac{e^{itu}-1}{u}\,du-\frac{e^{it\frac{k}{n-m}}-1}{(n-m)\frac{k}{n-m}}\Big|\\&=\Big|\int_{\frac{k-1-m}{n-m}}^{\frac{k-m}{n-m}} \eta_t(u)\,du-\frac{\eta_t\Big(\frac{k}{n-m}\Big)}{n-m}\Big|=\Big|\int_{\frac{k-1-m}{n-m}}^{\frac{k-m}{n-m}} \Big\{\eta_t(u)-\eta_t\Big(\frac{k}{n-m}\Big)\Big\}\,du\Big|\\&\leqslant \frac{1}{n-m}\sup_{u \in [\frac{k-1-m}{n-m},\frac{k-m}{n-m}]}\Big|\eta_t(u)-\eta_t\Big(\frac{k}{n-m}\Big)\Big|\\&\leqslant  \frac{1}{n-m}\cdot \Big\{\sup_{u \in [\frac{k-1-m}{n-m},\frac{k-m}{n-m}]}\Big|u-\frac{k}{n-m}\Big|\Big\}\cdot\Big\{\sup_{x \in \mathbb{R}}\big|\eta^\prime_t(x)\big|\Big\}=\frac{m+1}{(n-m)^2}\cdot\sup_{x \in \mathbb{R}}\big|\eta^\prime_t(x)\big|\\&\leqslant \frac{(m+1)t^2}{2(n-m)^2},\end{align*}
since
$$\sup_{x \in \mathbb{R}}\big|\eta^\prime_t(x)\big|\leqslant\frac{t^2}{2},$$
as we are going to prove. First
\begin{align*}
\big|\eta^\prime_t(x)\big|^2=\Big|\frac{it x e^{itx}- e^{itx}+1}{x^2}\Big|^2=\frac{\delta(tx)}{x^4},
\end{align*}
where $\delta(u) =2(1-u \sin u -\cos u)+u^2$. Put now $H(u) = \frac{u^4}{4}- \delta (u)$. We have
$$H^\prime(u) = u^3 - 2u(1-\cos u)= 4u\Big(\frac{u^2}{4} - \sin^2 \frac{u}{2}\Big)\geqslant 0, \qquad u\geqslant 0;$$
hence $H$ is non--decreasing for $u\geqslant 0$, and from the fact that  $H(0)=0$, we deduce that $H(u) \geqslant 0$ for every $u\geqslant 0$, hence also for every $u \in \mathbb{R}$ since $H(-u) =H(u)$. In other words $\delta(u) \leqslant\frac{u^4}{4}$ and as a consequence
\begin{align*}
\big|\eta^\prime_t(x)\big|^2=\frac{\delta(tx)}{x^4} \leqslant\frac{t^4}{4}.
\end{align*}

\hfill
 $\Box$
 
 \bigskip
 \noindent
The following result specifies Proposition 1 of \cite{HT} quantitatively in terms of distribution functions.
\begin{proposition}\label{W2}
There exists an absolute positive constant $C$ such that, for all positive integers $n$, $m$, with $n>m \geqslant 2$,
\begin{equation*}\sup_{x \in \mathbb{R}}\Big|P\Big(\frac{T_m^n}{n-m}\leqslant x\Big)- D(x)\Big|\leqslant C g_{m,n},\end{equation*} where \begin{equation*}g_{m,n}=\exp\Bigg(C\Big\{\frac{\log\frac{n}{m}}{(n-m)^2}+ \frac{m+2}{n-m}\Big\}\log^2\frac{n}{m}\Bigg)-1+ \frac{1}{\log\frac{n}{m}}.
\end{equation*}
\end{proposition}

\noindent
{\it Proof.}
In view of Theorem 2 p. 109 in \cite{P}, if $\tau$ is an arbitrary positive number, then for $b > \frac{1}{2 \pi}$

\begin{align*}
&\sup_{x \in \mathbb{R}}\Big|P\Big(\frac{T_m^n}{n-m}\leqslant x\Big)- D(x)\Big|\leqslant b \int_{-\tau}^\tau \frac{\big|\exp\Big\{\int_0^1 \frac{e^{itu}-1}{u}\,du\Big\}-\phi_{\frac{T_m^n}{n-m}}(t)\big|}{|t|}
\,{\rm d}t+\frac{r(b)}{\tau} \sup_{x \in \mathbb{R}}|D^\prime (x)|
\end{align*}
where $r(b)$ is a positive constant depending on $b$ only. Hence
\begin{align*}
&\sup_{x \in \mathbb{R}}\Big|P\Big(\frac{T_m^n}{n-m}\leqslant x\Big)- D(x)\Big|\\&\leqslant C\Big\{ \int_{-\tau}^\tau \frac{\big|\exp\Big\{\int_0^1 \frac{e^{itu}-1}{u}\,du\Big\}-\phi_{\frac{T_m^n}{n-m}}(t)\big|}{|t|}
\,{\rm d}t+\frac{1}{\tau} \sup_{x \in \mathbb{R}}|D^\prime (x)|\Big\}\\& \leqslant   C\Big\{ \int_{-\tau}^\tau\frac{f_{m,n}(t)}{|t|}
\,{\rm d}t+\frac{1}{\tau} \Big\},\end{align*}
by Proposition \ref{W1}.
Now, for every positive constant $A$  we have
\begin{align*}
&\sup_{0 \leqslant x \leqslant \tau }\frac{e^{A  x^2}-1}{x}
= \frac{e^{A  \tau^2}-1}{\tau}.
\end{align*}
Applying this with $A =C\Big(\frac{\log\frac{n}{m}}{(n-m)^2}+ \frac{m+2}{n-m}\Big)$ we obtain
\begin{equation*}
\int_{-\tau}^\tau \frac{f_{m,n}(t)}{|t|}
\,{\rm d}t \leqslant 2\tau\frac{e^{C\Big(\frac{\log\frac{n}{m}}{(n-m)^2}+ \frac{m+2}{n-m}\Big) \tau^2}-1}{\tau}=2\Big(e^{C\Big(\frac{\log\frac{n}{m}}{(n-m)^2}+ \frac{m+2}{n-m}\Big)\tau^2}-1\Big).
\end{equation*}
Hence, for every $\tau$,
\begin{equation}\label{questa}
\sup_{x \in \mathbb{R}}\Big|P\Big(\frac{T_m^n}{n-m}\leqslant x\Big)- D(x)\Big|\leqslant C\Bigg\{\exp\Bigg(C\Big\{\frac{\log\frac{n}{m}}{(n-m)^2}+ \frac{m+2}{n-m}\Big\}  \tau^2\Bigg)-1+\frac{1}{\tau}\Bigg\}.
\end{equation}
Taking $\tau= \log \frac{n}{m}$
we get
\begin{align*}
&\exp\Bigg(C\Big\{\frac{\log\frac{n}{m}}{(n-m)^2}+ \frac{m+2}{n-m}\Big\}   \tau^2\Bigg)-1+\frac{1}{\tau}\\ &=\exp\Bigg(C\Big\{\frac{\log\frac{n}{m}}{(n-m)^2}+ \frac{m+2}{n-m}\Big\}\log^2 \frac{n}{m}\Bigg)-1+ \frac{1}{\log\frac{n}{m}}=g_{m,n},
\end{align*}
so that, from \eqref{questa}
\begin{equation*}
\sup_{x \in \mathbb{R}}\Big|P\Big(\frac{T_m^n}{n-m}\leqslant x\Big)- D(x)\Big|\leqslant Cg_{m,n},
\end{equation*}
as claimed.

\hfill
 $\Box$
 
\section{Preliminaries on the Dickman distribution II: a convergence result}
This section is devoted to a convergence result for the characteristic functions of $\frac{T_n}{n}$ that has been used before in the proof of Proposition \ref{stimabase}; it gives also a weak form of the Local Limit Theorem (see Remark \ref{weak}).
\begin{proposition}\label{ZS} We have

\begin{itemize}
\item[\rm(a)]$$\int_{-\infty}^{+\infty}\big|\phi(t)\big|^2 \,dt < + \infty;$$

\item[\rm(b)]

\begin{equation*}\int_{-\pi n}^{\pi n}\Big|\phi_{\frac{T_n}{n}}(u)\Big|^2\,du\to \int_{-\infty}^{\infty}|\phi(u)|^2\,du, \qquad n \to \infty.
\end{equation*}\end{itemize}
\end{proposition}

 \noindent
{\it Proof.}

\bigskip
\noindent

\begin{itemize}
\item[(a)]

By symmetry ($t \mapsto \big|\phi(t)\big|^2$ is an even function), it is sufficient to prove that $$\int_{0}^{+\infty}\big|\phi(t)\big|^2 \,dt < + \infty.$$
Theorem 2 p. 11 of \cite{P} assures that there exist positive constants $\delta$ and $C$ such that
$$|\phi(t)|\leqslant 1-C t^2, \qquad |t|< \delta.$$
This implies that
$$\int_0^\delta |\phi(t)|^2 \, dt \leqslant \int_0^\delta (1-C t^2)^2 \,dt =C .$$
Let's turn to $\int_\delta^{+\infty} |\phi(t)|^2 \, dt.$ First observe that
\begin{align}\label{modulo}
 \nonumber\big|\phi(t)\big|^2 &= \phi(t)\phi(-t)= \exp \Big\{\int_0^1 \frac{e^{itu}-1}{u}\,du\Big\} \cdot \exp \Big\{\int_0^1 \frac{e^{-itu}-1}{u}\,du\Big\}\\& = \exp \Big\{-2\int_0^1 \frac{1-\cos tu}{u}\,du\Big\}.
\end{align}Now, for every $\epsilon \in (0,t)$ 
\begin{align*}
&\int_0^1\frac{1-\cos tu}{u}du = \int_0^t\frac{1-\cos z}{z}dz\geqslant \int_\epsilon^t\frac{1-\cos z}{z}dz= \Big[\frac{z-\sin z}{z}\Big]_\epsilon^t + \int_\epsilon^t \frac{z-\sin z}{z^2}\, dz \\&= \log\frac{t}{\epsilon}- \frac{\sin t}{t}+\frac{\sin \epsilon}{\epsilon}- \int_\epsilon^t\frac{\sin z}{z^2}\, dz\geqslant\log\frac{t}{\epsilon}+C
\end{align*}
(the constant $C$ might be negative here, but this is irrelevant as it appears clearly from the sequel).
Hence, taking $\epsilon = \delta$, 

\begin{align*}
&\int_\delta^{+\infty} |\phi(t)|^2 \, dt\leqslant \int_\delta^{+\infty}\exp\Big\{-2 \Big(\log\frac{t}{\delta}+C\Big)\Big\}\, dt\leqslant C\int_\delta^{+\infty}\frac{1}{t^2}\, dt = C.
\end{align*}

\item[(b)]

By part (a), the Proposition will be proved if we show that
\begin{equation}\label{asserto}
\int_{-\pi n}^{\pi n}\Big\{\Big|\phi_{\frac{T_n}{n}}(u)\Big|^2-|\phi(u)|^2\Big\}\,du\to 0, \qquad n \to \infty.
\end{equation}
Recall that, by  Proposition \ref{law}, $\Big|\phi_{\frac{T_n}{n}}\Big|$ converges to $|\phi|$ pointwise and uniformly on every bounded interval. Hence, for any positive $A$, 
 \begin{equation*}
\int_{-A}^{A}\Big\{\Big|\phi_{\frac{T_n}{n}}(u)\Big|^2-|\phi(u)|^2\Big\}\,du\to 0, \qquad n \to \infty.
\end{equation*}
Thus we are left with the proof of
\begin{equation}\label{uno}
\int_{\{A \leqslant |t|\leqslant n \pi\}}\Big\{\Big|\phi_{\frac{T_n}{n}}(u)\Big|^2-|\phi(u)|^2\Big\}\,du\to 0, \qquad n \to \infty.
\end{equation}
We split the first member of \eqref{uno} as follows: for a fixed $\epsilon\in (0,1)$, write
$$\int_{\{A \leqslant |t|\leqslant n \pi\}}= \int_{\{A \leqslant |t|\leqslant\epsilon\pi\sqrt[5]{n} \}}+\int_{\{\epsilon\pi\sqrt[5]{n}  \leqslant |t|\leqslant n \pi\}} = I_1 + I_2.$$
We consider the two summands $I_1$ and $I_2$ separately.\begin{itemize}
\item[($I_1$)] Notice that
\begin{align}\label{trigonometrica}&\nonumber\Big|\phi_{T_n}(t)\Big|^2= \prod_{k=1}^n \Big|1+ \frac{e^{ikt}-1}{k}\Big|^2 = \prod_{k=1}^n \Big\{\Big(1- \frac{1}{k}+\frac{1}{k} \cos kt\Big)^2+\Big(\frac{1}{k} \sin kt\Big)^2\Big\}\\&= \prod_{k=1}^n \Big\{1+ \frac{2(k-1)}{k^2}\Big( \cos kt-1\Big)\Big\}=\exp \Big\{\sum_{k=1}^n\log \Big[1+\frac{2(k-1)}{k^2}\Big(\cos kt-1\Big)\Big]\Big\}.\end{align}

Hence
\begin{align*}|I_1| \leqslant &\int_{\{A \leqslant |t|\leqslant\epsilon\pi\sqrt[5]{n}\}}\Bigg|\Big|\phi_{\frac{T_n}{n}}(t)\Big|^2-\exp \Big\{\sum_{k=1}^n\frac{2(k-1)}{k^2}\Big(\cos \frac{kt}{n} -1\Big)\Big\}\Bigg|\,dt+\\& +\int_{\{A \leqslant |t|\leqslant\epsilon\pi\sqrt[5]{n} \}}\Bigg|\exp \Big\{\sum_{k=1}^n\frac{2(k-1)}{k^2}\Big(\cos \frac{kt}{n} -1\Big)\Big\}-\exp \Big\{\sum_{k=1}^n\frac{2}{k}\Big(\cos \frac{kt}{n} -1\Big)\Big\}\Bigg|\,dt+\\& +\int_{\{A \leqslant |t|\leqslant\epsilon\pi\sqrt[5]{n} \}}\Bigg|\exp \Big\{\sum_{k=1}^n\frac{2}{k}\Big(\cos \frac{kt}{n} -1\Big)\Big\}-|\phi(t)|^2\Bigg|\,dt= I_{11}+ I_{12}+I_{13}.
\end{align*}
We consider the three summand $ I_{11}$, $ I_{12}$ and $I_{13}$ separately.
\begin{itemize}
\item[($I_{11}$)]First observe that, by relation \eqref{trigonometrica}, 
\begin{align*}
I_{11}&=\int_{\{A \leqslant |t|\leqslant\epsilon\pi\sqrt[5]{n} \}}\Bigg|\Big|\phi_{\frac{T_n}{n}}(t)\Big|^2-\exp \Big\{\sum_{k=1}^n\frac{2(k-1)}{k^2}\Big(\cos \frac{kt}{n} -1\Big)\Big\}\Bigg|\,dt\\ &=\int_{\{A \leqslant |t|\leqslant\epsilon\pi\sqrt[5]{n} \}}\Bigg|\exp \Big\{\sum_{k=1}^n\log \Big[1+\frac{2(k-1)}{k^2}\Big(\cos \frac{kt}{n}-1\Big)\Big]\Big\}\\ &-\exp \Big\{\sum_{k=1}^n\frac{2(k-1)}{k^2}\Big(\cos \frac{kt}{n} -1\Big)\Big\}\Bigg|\,dt\\& \leqslant 
\int_{\{A \leqslant |t|\leqslant\epsilon\pi\sqrt[5]{n}\}}\Bigg|\exp \Bigg\{\sum_{k=1}^n\Big\{\log \Big[1+\frac{2(k-1)}{k^2}\Big(\cos \frac{kt}{n}-1\Big)\Big]\\ &- \frac{2(k-1)}{k^2}\Big(\cos \frac{kt}{n} -1\Big)\Big\}\Bigg\}-1\Bigg|\,dt,\end{align*}
 since
$$0 \leqslant \exp \Big\{\sum_{k=1}^n\frac{2(k-1)}{k^2}\Big(\cos \frac{kt}{n} -1\Big)\Big\}\leqslant 1.$$
Now using the development
$$\log(1+z) - z = \sum_{j \geqslant 2}\frac{(-1)^j}{j}z^j,\quad |z|< 1$$
with $z=\frac{2(k-1)}{k^2}\big(\cos \frac{kt}{n} -1\big)$ (which, for sufficiently large $n$, is 
strictly less than 1 in modulus for every $k \geqslant 1$ ) we get
\begin{align*}
&\log \Big[1+\frac{2(k-1)}{k^2}\Big(\cos \frac{kt}{n}-1\Big)\Big]-\frac{2(k-1)}{k^2}\Big(\cos \frac{kt}{n} -1\Big)\\&=
\sum_{j \geqslant 2}\frac{(-1)^j2^j(k-1)^j}{j\cdot k^{2j}}\Big(\cos \frac{kt}{n} -1\Big)^j=\sum_{j \geqslant 2}\frac{2^j(k-1)^j}{j k^{2j}}\Big(1-\cos \frac{kt}{n} \Big)^j\geqslant 0.\end{align*}
It follows that
\begin{align*}&\int_{\{A \leqslant |t|\leqslant\epsilon\pi\sqrt[5]{n} \}}\Bigg|\exp \left\{\sum_{k=1}^n\Big\{\log \Big[1+\frac{2(k-1)}{k^2}\Big(\cos \frac{kt}{n}-1\Big)\Big]-\frac{2(k-1)}{k^2}\Big(\cos \frac{kt}{n} -1\Big)\Big\}\right\}-1\Bigg|\,dt\\&=\int_{\{A \leqslant |t|\leqslant\epsilon\pi\sqrt[5]{n} \}}\Bigg(\exp \Big\{\sum_{k=1}^n\sum_{j \geqslant 2}\frac{2^j(k-1)^j}{j k^{2j}}\Big(1-\cos \frac{kt}{n} \Big)^j\Big\}-1\Bigg)\,dt.\end{align*}
Using the inequality $1-\cos z \leqslant z^2$ the above can be bounded by
\begin{align*}
&\int_{\{A \leqslant |t|\leqslant\epsilon\pi\sqrt[5]{n} \}}\Bigg(\exp \Big\{\sum_{k=1}^n\sum_{j \geqslant 2}\frac{2^j(k-1)^j}{j k^{2j}}\Big( \frac{kt}{n} \Big)^{2j}\Big\}-1\Bigg)\,dt\\&=\int_{\{A \leqslant |t|\leqslant\epsilon\pi\sqrt[5]{n} \}}\Bigg(\exp \Big\{\sum_{k=1}^n\sum_{j \geqslant 2}\frac{2^j(k-1)^j}{j }\Big( \frac{t}{n} \Big)^{2j}\Big\}-1\Bigg)\,dt
\\&=\int_{\{A \leqslant |t|\leqslant\epsilon\pi\sqrt[5]{n} \}}\Bigg(\exp \Big\{\sum_{j \geqslant 2}\frac{1}{j }\Big( \frac{\sqrt 2t}{n} \Big)^{2j}\sum_{k=1}^n(k-1)^j\Big\}-1\Bigg)\,dt\\&\leqslant
\int_{\{A \leqslant |t|\leqslant\epsilon\pi\sqrt[5]{n} \}}\Bigg(\exp \Big\{\sum_{j \geqslant 2}\frac{1}{j }\Big( \frac{\sqrt 2t}{n} \Big)^{2j}\big(\int_{1}^nx^j\,dx\big)\Big\}-1\Bigg)\,dt\\&\leqslant \int_{\{A \leqslant |t|\leqslant\epsilon\pi\sqrt[5]{n} \}}\Bigg(\exp \Big\{\sum_{j \geqslant 2}\frac{1}{j }\Big( \frac{\sqrt 2t}{n} \Big)^{2j}\Big(\frac{n^{j+1}}{j+1}\Big)\Big\}-1\Bigg)\,dt\\&=\int_{\{A \leqslant |t|\leqslant\epsilon\pi\sqrt[5]{n} \}}\Bigg(\exp \Big\{\sum_{j \geqslant 2}\frac{1}{j(j+1) }\cdot \frac{(\sqrt 2t)^{2j}}{n^{j-1}} \Big\}-1\Bigg)\,dt\\&=\int_{\{A \leqslant |t|\leqslant\epsilon\pi\sqrt[5]{n} \}}\Bigg(\exp \Big\{C\frac{t^4}{n}\sum_{j \geqslant 0}\frac{1}{(j+3)(j+2) }\cdot \frac{(\sqrt 2t)^{2j}}{n^{j}} \Big\}-1\Bigg)\,dt.\end{align*}
Now, for $|t|\leqslant\epsilon\pi\sqrt[5]{n}$ we have also$|t|\leqslant\pi\sqrt[5]{n}$ (recall that $\epsilon <1$); hence there exists $n_0$ not dependent on $\epsilon$ such that, for $n>n_0$
$$\frac{(\sqrt 2t)^{2}}{n}\leqslant \frac{C}{n^{\frac{3}{5}}}\leqslant \frac{1}{2}.$$
 Hence, for $n>n_0$,
$$\sum_{j \geqslant 0}\frac{1}{(j+3)(j+2) }\cdot \frac{(\sqrt 2t)^{2j}}{n^{j}}\leqslant \sum_{j \geqslant 0}\frac{1}{(j+3)(j+2) }\cdot \frac{1}{2^j} =C,$$
and we get
\begin{align*}
&\int_{\{A \leqslant |t|\leqslant\epsilon\pi\sqrt[5]{n} \}}\Bigg(\exp \Big\{C\frac{t^4}{n}\sum_{j \geqslant 0}\frac{1}{(j+3)(j+2)\cdot }\cdot \frac{(\sqrt 2t)^{2j}}{n^{j}} \Big\}-1\Bigg)\,dt \\&\leqslant \int_{\{A \leqslant |t|\leqslant\epsilon\pi\sqrt[5]{n} \}}\Bigg(\exp \Big\{C\frac{t^4}{n}\Big\}-1\Bigg)\,dt\leqslant \Bigg(\exp \Big\{C\frac{(\epsilon\pi\sqrt[5]{n} )^4}{n}\Big\}-1\Bigg)\epsilon\pi\sqrt[5]{n}\\&=C \cdot \frac{e^ {\frac{C}{\sqrt[5]{n}}}-1}{\frac{C}{\sqrt[5]{n}}}\epsilon^5< C \epsilon^5.
\end{align*}
\item[($I_{12}$)]Here we observe that$$0 \leqslant \exp \Big\{\sum_{k=1}^n\frac{2}{k}\Big(\cos \frac{kt}{n} -1\Big)\Big\}\leqslant 1,$$
hence
\begin{align*}
&I_{12}=\int_{\{A \leqslant |t|\leqslant\epsilon\pi\sqrt[5]{n} \}}\Bigg|\exp \Big\{\sum_{k=1}^n\frac{2(k-1)}{k^2}\Big(\cos \frac{kt}{n} -1\Big)\Big\}-\exp \Big\{\sum_{k=1}^n\frac{2}{k}\Big(\cos \frac{kt}{n} -1\Big)\Big\}\Bigg|\,dt\\&\leqslant\int_{\{A \leqslant |t|\leqslant\epsilon\pi\sqrt[5]{n} \}}\Bigg|\exp \Big\{\sum_{k=1}^n\Big(\frac{2(k-1)}{k^2}- \frac{2}{k}\Big)\Big(\cos \frac{kt}{n} -1\Big)\Big\}-1\Bigg|\,dt\\&=\int_{\{A \leqslant |t|\leqslant\epsilon\pi\sqrt[5]{n} \}}\Bigg(\exp \Big\{\sum_{k=1}^n\frac{2}{k^2}\Big(1-\cos \frac{kt}{n} \Big)\Big\}-1\Bigg)\,dt\\& \leqslant \int_{\{A \leqslant |t|\leqslant\epsilon\pi\sqrt[5]{n} \}}\Bigg(\exp \Big\{\sum_{k=1}^n\frac{2}{k^2}\Big( \frac{kt}{n} \Big)^2\Big\}-1\Bigg)\,dt = \int_{\{A \leqslant |t|\leqslant\epsilon\pi\sqrt[5]{n} \}}\Big(e^{\frac{2t^2}{n}}-1\Big)\,dt\\& \leqslant \Big(e^{\frac{2(\epsilon\pi\sqrt[5]{n})^2}{n}}-1\Big)\epsilon\pi\sqrt[5]{n}= \frac{e^{\frac{2(\epsilon\pi)^2}{n^{3/5}}}-1}{\frac{2(\epsilon\pi)^2}{n^{3/5
 }}}\cdot \frac{2(\epsilon\pi)^2}{n^{3/5}}\epsilon\pi\sqrt[5]{n}\to 0 , \qquad n \to \infty.
\end{align*}
\item[($I_{13}$)] Recalling the explicit form of $|\phi(t)|^2$ given in equation \eqref{modulo}, we have
\begin{align*}
&I_{13}=\int_{\{A \leqslant |t|\leqslant\epsilon\pi\sqrt[5]{n} \}}\Bigg|\exp \Big\{\sum_{k=1}^n\frac{2}{k}\Big(\cos \frac{kt}{n} -1\Big)\Big\}-|\phi(t)|^2\Bigg|\,dt \\&\leqslant  \int_{\{A \leqslant |t|\leqslant\epsilon\pi\sqrt[5]{n} \}}\Bigg|\exp \Big\{\sum_{k=1}^n\frac{2}{k}\Big(\cos \frac{kt}{n} -1\Big)\Big\}-\exp \Big\{2\int_{0}^1\frac{\cos tu -1}{u}\,du\Big\}\Bigg|\,dt,\end{align*}
and, observing again that
$$0\leqslant \exp \Big\{\sum_{k=1}^n\frac{2}{k}\Big(\cos \frac{kt}{n} -1\Big)\Big\}\leqslant 1,$$
we get
\begin{align*}
 &I_{13}\leqslant \int_{\{A \leqslant |t|\leqslant\epsilon\pi\sqrt[5]{n} \}}\Bigg|\exp \Big\{2\int_{0}^1\frac{\cos tu -1}{u}\,du-2\sum_{k=1}^n\frac{n}{k}\Big(\cos \frac{kt}{n} -1\Big)\frac{1}{n}\Big\}-1\Bigg|\,dt\\&= \int_{\{A \leqslant |t|\leqslant\epsilon\pi\sqrt[5]{n} \}}\Bigg|\exp \Big\{2\sum_{k=1}^n\int_{\frac{k-1}{n}}^{\frac{k}{n}}\big[\gamma_t(u)-\gamma_t\Big(\frac{k}{n}\Big)\big]\,du\Big\}-1\Bigg|\,dt\\&\leqslant \int_{\{A \leqslant |t|\leqslant\epsilon\pi\sqrt[5]{n} \}}\Bigg(\exp \Big\{2\sum_{k=1}^n\int_{\frac{k-1}{n}}^{\frac{k}{n}}\big|\gamma_t\Big(\frac{k}{n}\Big)-\gamma_t(u)\big|\,du\Big\}-1\Bigg)\,dt,
\end{align*}
where we have put
$$\gamma_t(u)=\frac{\cos tu -1}{u}$$
and have used the inequality $|e^x - 1|\leqslant e^{|x|}-1.$

It is not difficult to see that $$\sup_{u \in \mathbb{R}}\big|\gamma^\prime_t(u)\big|= C t^2;$$
in fact $\gamma^\prime_t(u)= \eta(tu) t^2$, with
$$\eta(z)= \frac{1-z \sin z -\cos z }{z^2}, $$
and proving that $\sup_{z \in \mathbb{R}}|\eta(z)|= C < +\infty$
is a simple exercise. Thus, by Lagrange's Theorem,
\begin{equation}\label{boundintegrale}
\int_{\frac{k-1}{n}}^{\frac{k}{n}}\big|\gamma_t\Big(\frac{k}{n}\Big)-\gamma_t(u)\big|du\leqslant C t^2 \int_{\frac{k-1}{n}}^{\frac{k}{n}}\Big|\frac{k}{n}-u\Big|du\leqslant C   \frac{t^2}{n^2}.
\end{equation}
Using \eqref{boundintegrale} in the last bound for $I_{13}$ we find
\begin{align*}
&I_{13}\leqslant \int_{\{A \leqslant |t|\leqslant\epsilon\pi\sqrt[5]{n} \}}\Bigg(\exp \Big\{2\sum_{k=1}^n\int_{\frac{k-1}{n}}^{\frac{k}{n}}\big|\gamma_t\Big(\frac{k}{n}\Big)-\gamma_t(u)\big|\,du\Big\}-1\Bigg)\,dt\\&
\leqslant \int_{\{A \leqslant |t|\leqslant\epsilon\pi\sqrt[5]{n} \}}\Bigg(\exp \Big\{C \frac{t^2}{n}\Big\}-1\Bigg)\,dt\leqslant \Bigg(\exp \Big\{C \frac{(\epsilon\pi\sqrt[5]{n})^2}{n}\Big\}-1\Bigg)\epsilon\pi\sqrt[5]{n}\\& =C
\frac{e^{\frac{C}{n^{3/5}}}-1}{ \frac{C}{n^{3/5}}}\frac{C}{n^{2/5}} \to 0, \qquad n \to \infty.
\end{align*}

\end{itemize}

\item[($I_2$)]We recall that
\begin{align*}
I_2&=\int_{\{\epsilon\pi\sqrt[5]{n}  \leqslant |t|\leqslant n \pi\}}
\Big\{\Big|\phi_{\frac{T_n}{n}}(u)\Big|^2-|\phi(u)|^2\Big\}\,du\\& \leqslant \int_{\{\epsilon\pi\sqrt[5]{n}  \leqslant |t|\leqslant n \pi\}}
\Big|\phi_{\frac{T_n}{n}}(u)\Big|^2\,du+\int_{\{\epsilon\pi\sqrt[5]{n}  \leqslant |t|\leqslant n \pi\}}
|\phi(u)|^2\,du.\end{align*}

The second summand above goes to 0 as $n \to\infty$ since $|\phi(u)|^2$ is integrable on $\mathbb{R}$ (recall point (a) of this proposition); hence we have to prove that
$$\int_{\{\epsilon\pi\sqrt[5]{n}  \leqslant |t|\leqslant n \pi\}}
\Big|\phi_{\frac{T_n}{n}}(u)\Big|^2\,du\to 0, \qquad n \to \infty.$$
By relation \eqref{trigonometrica}, we have, for every $k_0 \in \mathbb{N}$ and $n \geqslant k_0 $
\begin{align*}
&\Big|\phi_{T_n}(u)\Big|^2\\& = \exp \Big\{\sum_{k=1}^{k_0-1 }\log \Big[1+\frac{2(k-1)}{k^2}\Big(\cos kt-1\Big)\Big]\Big\}\cdot\exp \Big\{\sum_{k=k_0}^{n }\log \Big[1+\frac{2(k-1)}{k^2}\Big(\cos kt-1\Big)\Big]\Big\}\\&\leqslant
\exp \Big\{\sum_{k=k_0}^{n }\log \Big[1+\frac{2(k-1)}{k^2}\Big(\cos kt-1\Big)\Big]\Big\}.\end{align*}
Now, using the relation
$$\big|\log(1-z)+z\big|\leqslant |z|^2, \qquad |z|< \frac{1}{2}, \qquad z \in \mathbb{C}$$
with $z= -\frac{2(k-1)}{k^2}\big(\cos kt-1\big)$ and choosing $k_0$ such that $\big| \frac{2(k-1)}{k^2}\big(\cos kt-1\big)\big|< \frac{1}{2}$ for $k> k_0$ and every $t$, we find
\begin{align*}
\log \Big[1+\frac{2(k-1)}{k^2}\Big(\cos kt-1\Big)\Big]\leqslant -
\frac{2(k-1)}{k^2}\Big(1-\cos kt\Big)+\frac{4(k-1)^2}{k^4}\Big(1-\cos kt\Big)^2.\end{align*}
Thus, by the obvious inequality  $1-\cos x \leqslant 2$,
\begin{align*}
&\exp \Big\{\sum_{k=k_0}^{n }\log \Big[1+\frac{2(k-1)}{k^2}\Big(\cos kt-1\Big)\Big]\Big\}\\&\leqslant\exp \Big\{-\sum_{k=k_0}^{n }
\frac{2(k-1)}{k^2}\Big(1-\cos kt\Big)+\sum_{k=k_0}^{n }\frac{4(k-1)^2}{k^4}\Big(1-\cos kt\Big)^2\Big\}\\&\leqslant
\exp \Big\{-2\sum_{k=k_0}^{n }
\frac{k-1}{k^2}+2\sum_{k=k_0}^{n }
\frac{k-1}{k^2}\cos kt+\sum_{k=k_0}^{n }\frac{16(k-1)^2}{k^4}\Big\}\\&\leqslant\exp \Big\{-2\sum_{k=k_0}^{n }
\frac{1}{k}+2\sum_{k=k_0}^{n}
\frac{\cos kt}{k}+C\Big\}\\&\leqslant\exp \Big\{-2\int_{k_0}^n\frac{1}{x}\, dx+2\sum_{k=k_0}^{n}
\frac{\cos kt}{k}+C\Big\}= \frac{C}{n^2} \exp \Big\{2\sum_{k=k_0}^{n}
\frac{\cos kt}{k}\Big\}.\end{align*}
By \cite{Z}, p. 191 we know that
$$\sup_{n \geqslant 1}\Big|\sum_{k=1}^{n}
\frac{\cos kt}{k}\Big|\leqslant \log \frac{1}{t}+C.$$
Hence
$$\sup_{n \geqslant k_0}\Big|\sum_{k=k_0}^{n}
\frac{\cos kt}{k}\Big|\leqslant  \sup_{n \geqslant 1}\Big|\sum_{k=1}^{n}
\frac{\cos kt}{k}\Big|+ \Big|\sum_{k=1}^{k_0-1}
\frac{\cos kt}{k}\Big|\leqslant \log \frac{1}{t}+C.$$
 It follows that
$$ \frac{C}{n^2} \exp \Big\{2\sum_{k=k_0}^{n}
\frac{\cos kt}{k}\Big\}\leqslant \frac{C}{n^2}\exp \Big\{2\Big(\log \frac{1}{t}+C\Big)\Big\}= \frac{C}{n^2t^2},$$
so that
\begin{align*}
&\int_{\{\epsilon\pi\sqrt[5]{n}  \leqslant |t|\leqslant n \pi\}}
\Big|\phi_{\frac{T_n}{n}}(u)\Big|^2\,du= n \int_{\frac{\epsilon}{n^{4/5}}}^{\pi}\Big|\phi_{T_n}(t)\Big|^2\,dt\leqslant n \int_{\frac{\epsilon}{n^{4/5}}}^{\pi}\frac{C}{n^2t^2}\,dt\\&\leqslant 
C \frac{n^{4/5}}{n}= \frac{C}{n^{1/5}}\to 0, \qquad n \to \infty.\end{align*}

\end{itemize}

\end{itemize}

Now \eqref{asserto} is proved by letting $\epsilon \to 0$ in the estimation of $I_{11}$.

\hfill $\Box$

\begin{remark}\label{weak}\rm 
The relation \eqref{asserto} yields a weak form of  Local Limit Theorem (see Corollary \ref{local}). Let $\widehat{T_n}$ be the $n$-th partial sum $$\widehat{T_n}=\sum_{k = 1}^n k \widehat{Z_k},$$ where $(\widehat{Z_n})_{n\geqslant 1}$ is an independent copy of $(Z_n)_{n\geqslant 1}$. Denote by $d_s$ the symmetrized Dickman density, which has characteristic function $\big|\phi\big|^2$. Then, by the inversion formula, \eqref{asserto} and Proposition \ref{ZS},

\begin{align*}
&2 \pi\big\{ n P(T_n -\widehat{T_n}= \kappa_n)- d_s(n^{-1}\kappa_n)\big\}\\&= \int_{-n \pi}^{n\pi}e^{-itn^{-1}\kappa_n }\Big\{\Big|\phi_{\frac{T_n}{n}}(t)\Big|^2-|\phi(t)|^2\Big\}\,dt+\int_{\{|t|> n\pi\}}e^{-itn^{-1}\kappa_n }|\phi(t)|^2
\,dt\to 0, \quad n \to \infty.\end{align*}

\end{remark}

\section{The correlation inequality}
In this section we present a correlation inequality for the sequence of random variables $(Y_n)_{n \geqslant 1}$, where
\begin{equation}\label{definizione}
Y_n = n 1_{\{T_n = \kappa_n\}}.
\end{equation}
\begin{theorem} {\bf(Basic correlation inequality).}\label{covarianzabase} Let $x>0$ be given and  let ${\boldsymbol \kappa }=(\kappa_n)$ be any fixed sequence of integers with $\lim_{n \to \infty} \frac{\kappa_n}{n}=x$. 
Then, for every $x>0$ and for $n> m \geqslant 2$
\begin{align*}
&|Cov(Y_m, Y_n)|\leqslant C \bigg\{\frac{n}{n-m}\chi^{({\boldsymbol \kappa },x)}_{m,n}+\frac{m}{n-m}+\chi^{({\boldsymbol \kappa },x)}_{2,n}+\frac{1}{n}\bigg\}.
\end{align*}  
where $C$ is a positive constant (depending on $x$) and
$$\chi_{m,n}^{({\boldsymbol \kappa },x)}=\frac{n-m}{\kappa_n-\kappa_m}\cdot\frac{\log\frac{n}{m}}{\sqrt {n-m}}+ \frac{n-m}{\kappa_n-\kappa_m}\cdot g_{m,n} +x\Big|\frac{n-m}{\kappa_n-\kappa_m}-\frac{1}{x}\Big|+ \frac{m+1}{\kappa_n-\kappa_m}.$$
\end{theorem}
{\it Proof.}
\begin{align*}
&Cov(Y_m, Y_n)= nm\big\{ P(T_m = \kappa_m,T_n = \kappa_n)-P(T_m = \kappa_m)P(T_n = \kappa_n)\big\}\\&= nm\big\{ P(T_m = \kappa_m, T_{m}^n= \kappa_n -\kappa_m)-P(T_m = \kappa_m)P(T_n = \kappa_n)\big\}\\&=\big\{mP(T_m = \kappa_m)\big\}\big\{n P(T_{m}^n= \kappa_n -\kappa_m)-nP(T_n = \kappa_n)\big\}.
\end{align*}
Hence, by the local Theorem (Corollary \ref{local}),  we have
\begin{align} \label{prima}
&\nonumber \big|Cov(Y_m, Y_n)\big| \leqslant C \big|n P(T_{m}^n= \kappa_n -\kappa_m)-nP(T_n = \kappa_n)\big|\\&\nonumber \leqslant C\Big(\Big|\frac{n}{n-m}   \big\{(n-m)P(T_{m}^n= \kappa_n -\kappa_m)-D^\prime(x)\big\}\Big|+\Big|\Big(\frac{n}{n-m}  -1\Big)D^\prime(x)\Big|\\&\nonumber +\big|nP(T_n = \kappa_n)-D^\prime (x)\big|\Big) \\&\nonumber \leqslant  C\Big(\Big|\frac{n}{n-m}   \big\{(n-m)P(T_{m}^n= \kappa_n -\kappa_m)-D^\prime(x)\big\}\Big|+\frac{m}{n-m}  +\big|nP(T_n = \kappa_n)-D^\prime(x)\big|\Big)\\&= C\Big(\frac{n}{n-m}  \Gamma+\frac{m}{n-m}  +\Delta\Big),
\end{align}
where we have put for simplicity
$$\Gamma=|(n-m)P(T_{m}^n= \kappa_n -\kappa_m)-D^\prime(x)|, \qquad \Delta =|nP(T_n = \kappa_n)-D^\prime(x)\big|.$$
The aim is to obtain bounds for $\Gamma$ and $\Delta$.
 \begin{itemize}
\item[\bf (a) $\boldsymbol \Gamma.$] Set $$\Phi =P\big((\kappa_n- \kappa_m)-n \leqslant T_m^n \leqslant (\kappa_n- \kappa_m)-(m+1)\big);$$
$$\Psi =\Big| (\kappa_n -\kappa_m) P(T_n=\kappa_n -\kappa_m)-\Phi
 \Big|. $$
\begin{align*}
&\Gamma\leqslant\frac{n-m}{\kappa_n -\kappa_m}\Psi+\Big|\frac{n-m}{\kappa_n -\kappa_m}\Phi- D^\prime(x)\Big|=\frac{n-m}{\kappa_n -\kappa_m}\Psi+\Big|\frac{n-m}{\kappa_n -\kappa_m}\Phi- \frac{D(x)-D(x-1)}{x}\Big|,\end{align*}
by \eqref{equazione}.
From Proposition \ref{stimabase} we know that 
\begin{equation}\label{seconda}
\Psi \leqslant C\frac{\log \frac{n}{m}}{\sqrt {n-m}}.
\end{equation}
 Moreover, putting
 \begin{align*}
 &\Lambda=\Big|P\Big( \frac {T_m^n}{n-m} \leqslant \frac{(\kappa_n- \kappa_m)-(m+1)}{n-m}\Big)-D\Big(\frac{(\kappa_n- \kappa_m)-(m+1)}{n-m}\Big)\Big|,\\&\Theta =\Big|P\Big( \frac {T_m^n}{n-m}  \leqslant\frac{(\kappa_n- \kappa_m)-(n+1)}{n-m}\Big)-D\Big(\frac{(\kappa_n- \kappa_m)-(n+1)}{n-m}\Big)\Big|,\\&\Sigma =\Big|\frac{n-m}{\kappa_n-\kappa_m}D\Big(\frac{(\kappa_n- \kappa_m)-(m+1)}{n-m}\Big)-\frac{D(x)}{x}\Big|,\\&\Omega=
\Big|\frac{n-m}{\kappa_n-\kappa_m}D\Big(\frac{(\kappa_n-\kappa_m) -(n+1)}{n-m}\Big)-\frac{D(x-1)}{x}\Big|,
 \end{align*}
it is easily checked that
 \begin{align*}
  \Big|\frac{n-m}{\kappa_n -\kappa_m}\Phi- \frac{D(x)-D(x-1)}{x}\Big|\leqslant \frac{n-m}{\kappa_n -\kappa_m}\big(\Lambda+\Theta\big)+\Sigma+\Omega.
\end{align*}
We know from Proposition \ref{W2} that
$$\sup_{x \in \mathbb{R}}\Big|P\Big( \frac {T_m^n}{n-m} \leqslant x\Big)-D (x)\Big|\leqslant Cg_{m,n}.$$
Hence
\begin{align}\label{terza}\Lambda+\Theta
\leqslant Cg_{m,n} .
\end{align}
Moreover
\begin{align*}
&\Sigma\leqslant \frac{n-m}{\kappa_n-\kappa_m}\Big|D\Big(\frac{(\kappa_n- \kappa_m)-(m+1)}{n-m}\Big)-D(x)\Big| + \Big|\frac{n-m}{\kappa_n-\kappa_m}-\frac{1}{x}\Big|D(x) \\&\leqslant\frac{n-m}{\kappa_n-\kappa_m} \Big|D\Big(\frac{(\kappa_n- \kappa_m)-(m+1)}{n-m}\Big)-D(x)\Big| +\Big|\frac{n-m}{\kappa_n-\kappa_m}-\frac{1}{x}\Big|,
\end{align*}
and, by Lagrange Theorem, there exists $\xi_n$ such that
\begin{align}\label{quarta}
&\nonumber \frac{n-m}{\kappa_n-\kappa_m} \Big|D\Big(\frac{(\kappa_n- \kappa_m)-(m+1)}{n-m}\Big)-D(x)\Big| \leqslant\frac{n-m}{\kappa_n-\kappa_m}\Big|\frac{(\kappa_n- \kappa_m)-(m+1)}{n-m}-x\Big|D^\prime(\xi_n)
\\&\leqslant\frac{n-m}{\kappa_n-\kappa_m}\Big|\frac{(\kappa_n- \kappa_m)-(m+1)}{n-m}-x\Big|\leqslant x\Big|\frac{n-m}{\kappa_n-\kappa_m}-\frac{1}{x}\Big|+ \frac{m+1}{\kappa_n-\kappa_m},\end{align}
since $\sup_{x>0} D^\prime(x)=1$. For $\Omega$
we get exactly the same bound as in \eqref{quarta}.

\bigskip
\noindent
 In conclusion, from \eqref{seconda}, \eqref{terza} and \eqref{quarta} we have obtained
 \begin{align}\label{bound1}
 \Gamma=\big|(n-m)P(T_{m}^n= \kappa_n -\kappa_m)-D^\prime(x)\big|\leqslant C\chi^{({\boldsymbol \kappa },x)}_{m,n}.
 \end{align}
 \item[\bf (b) $\boldsymbol \Delta.$]
 Recall that
 $$\Delta=\big|nP(T_n = \kappa_n)-D^\prime(x)\big|.$$
 Notice that we cannot apply \eqref{bound1} directly since we have proved it for $m\geqslant 2$ only. Nevertheless, with $U= Z_1 + 2Z_2$,
\begin{align}\label{quinta}
&\nonumber \Delta=  \Big|\sum_{j=0}^3P(U=j)\Big(nP(T_2^n= \kappa_n -j)- D^\prime(x)\Big)\Big|\leqslant \sup_{0\leqslant j \leqslant 3}\Big|nP(T_2^n= \kappa_n -j)- D^\prime(x)\Big|
\\&\nonumber=\sup_{0\leqslant j \leqslant 3}\Big|\frac{n}{n-2}\big\{(n-2)P(T_2^n= \kappa_n -j)- D^\prime(x)\big\}+ \frac{2}{n-2}D^\prime(x)\Big|\\&\nonumber \leqslant\frac{n}{n-2}\Big\{\sup_{0\leqslant j \leqslant 3}|(n-2)P(T_2^n= \kappa_n -j)- D^\prime(x)|\Big\}+\frac{2}{n-2}\leqslant \frac{Cn}{n-2}\sup_{0\leqslant j \leqslant 3} \chi^{({\boldsymbol \kappa^{(j) },x)}}_{2,n}+\frac{2}{n-2}\\&\leqslant C \Big(\chi^{({\boldsymbol \kappa },x)}_{2,n}+\frac{1}{n}\Big),\end{align}
applying \eqref{bound1} (with $m=2$) for the sequence ${\boldsymbol \kappa^{(j)}}=(\kappa^{(j)}_n)_n$ defined as $\kappa^{(j)}_n = \kappa_n -j$ and noticing that  $\chi^{({\boldsymbol \kappa^{(j) },x)}}_{2,n}=\chi^{({\boldsymbol \kappa,x)}}_{2,n}$.
\end{itemize}

The two relations \eqref{bound1} and \eqref{quinta}, inserted into \eqref{prima}, conclude the proof.

\bigskip
  \hfill $\Box$

  \section{A general form of the Almost Sure Limit Theorem}

As we pointed out in the Introduction, the Almost Sure Limit Theorem that we are going to prove in the present section (i.e. Theorem \ref{principale}) is in the spirit of Theorem 1 of T. Mori's paper \cite{M}; in section 6 it will be applied to the sequence $(Y_n)$ defined in \eqref{definizione}: notice that Mori's result is not  applicable in the context of \eqref{definizione}, due to the fact that it requires that $\big|Cov(Y_m,Y_n)\big|\leqslant h\big(\frac{n}{m}\big)$  \underline{for all} $1 \leqslant m \leqslant n$ (for a suitable function $h$);  for $m=n$ this inequality  becomes $Var Y_m \leqslant h(1) =C$, i.e. the sequence $(Var Y_m)_{m \geqslant 1}$ must be bounded; unfortunately this is not true in our setting (see Lemma \ref{lemmaausiliario}).
\begin{theorem}\label{primoriusltato}
Let $(U_n)_{n \geqslant 1}$ a sequence of centered random variables. 
Assume that there exist two numbers $\alpha \geqslant 0$ and $\sigma >1$, a non--negative function  $f(u,z)$ defined on the set $\{u \geqslant 1, z \geqslant \sigma \}$, a non--negative double--indexed sequence $g$ defined on the set $\{(m,n)\in \mathbb{N}^2: \sigma m\leqslant n\}$ such that
$$
\sup_{n \geqslant \sigma m} g(m,n) = C < + \infty; \leqno \it(i)$$
(ii) 
uniformly in $u>0$ the functions $v \mapsto f\big(u, \frac{v}{u}\big)$ are ultimately non--increasing  (i.e.  there exists $m_0$ such that $v \mapsto f\big(u, \frac{v}{u}\big)$ is non--increasing on $(m_0,+\infty)$ and for every $u>0$);

\medskip
\noindent
(iii) the functions $$z \mapsto \phi(z) = \sup_{u \geqslant 1} \frac{f(u,z)}{z},\qquad u \mapsto F(u)=\int_\sigma^ u\phi(z) \, dz $$
are defined on $[\sigma, + \infty)$;

$$
|Cov(U_m, U_n)|\leqslant C\begin{cases}
m & \hbox{\rm for } m=n\\1& \hbox{\rm for } m<n \leqslant \sigma m;
\end{cases}\leqno \it(iv)
$$
(v) there exists $m_1$ such that, for $n> m \geqslant m_1$
$$
|Cov(U_m, U_n)|\leqslant g(m,n)\frac{1}{m^\alpha} f\Big(m,\frac{n}{m}\Big).
$$
Denote
$$V_n = \sum_{k=\sigma^{n-1}+1}^{\sigma^{n}} \frac{U_k}{k}.$$
Then, for every $n$ and every sufficiently large $m$
\begin{equation*}
E\big[(\sum_{i=m+1}^{m+n}V_i)^2\big]\leqslant C\left(n + \frac{1}{\sigma^{\alpha (m+n)}}\int_\sigma^{\sigma^n} F(u) u^{\alpha -1}\, du\right).
\end{equation*}
\end{theorem}

\bigskip
\noindent
{\it Proof.} Since
\begin{equation}\label{init}{\rm \bf E}\Big[\big(\sum_{i=m+1}^{m+n }V_i\big)^2\Big]=
\sum_{i=m+1}^{m+n}{\rm \bf E }[V_i^2]+2\sum_{m+1\leq i<j\leq
m+n}{\rm \bf E }[V_iV_j],\end{equation} we bound separately these
two summands. We have first \begin{align} \label{x}&\nonumber{\rm
\bf E }[V_i^2]= {\rm \bf E }\Big[\Big(\sum_{h= \sigma^{i-1}+1}^{\sigma^i
}\frac {U_h}{h}\Big) \Big(\sum_{k= \sigma^{i-1}+1}^{\sigma^i }\frac
{U_k}{k}\Big)\Big]=\sum_{h,k=\sigma^{i-1}+1}^{\sigma^i}\frac {1}{hk}{\rm \bf E
}[U_hU_k]\\&=\sum_{h=\sigma^{i-1}+1}^{\sigma^i}\frac {1}{h^2}{\rm \bf E
}[U_h^2]+2\sum_{\sigma^{i-1}+1\leq h<k\leq \sigma^i}\frac {1}{hk}{\rm \bf E
}[U_hU_k].
\end{align}
By the first inequality in (iv)
\begin{equation}\label{y}
\sum_{h=\sigma^{i-1}+1}^{\sigma^i}\frac {1}{h^2}{\rm \bf E
}[U_h^2] \leqslant C \sum_{h=\sigma^{i-1}+1}^{\sigma^i}\frac {1}{h}\leqslant C \log \frac{\sigma^i}{\sigma^{i-1}}= C.
\end{equation}
For the second summand in \eqref{x}, i.e.
$$\sum_{\sigma^{i-1}+1\leq h<k\leq \sigma^i}\frac {1}{hk}{\rm \bf E
}[U_hU_k],$$
we notice that
$$\frac{k}{\sigma} \leqslant \sigma^{i-1}<h$$
so that, by the second inequality in (iv), we have
\begin{equation}\label{z}
\sum_{\sigma^{i-1}+1\leq h<k\leq \sigma^i}\frac {1}{hk}{\rm \bf E
}[U_hU_k] \leqslant C\sum_{\sigma^{i-1}+1\leq h<k\leq \sigma^i}\frac {1}{hk}\leqslant C\Big(\sum_{k=\sigma^{i-1}+1}^{\sigma^{i}}\frac {1}{k}\Big)\Big(\sum_{h=\sigma^{i-1}+1}^{\sigma^{i}}\frac {1}{h}\Big)\leqslant C.
\end{equation}
The above relations \eqref{y} and \eqref{z}, used in \eqref{x}, give
\begin{equation}\label{primaparte}
\sum_{i=m+1}^{m+n}{\rm \bf E }[V_i^2]\leqslant  C\sum_{i=m+1}^{m+n}1=C n.
\end{equation}
Now we consider the second sum in \eqref{init}, i.e. $$\sum_{m+1\leq i<j\leq
m+n}{\rm \bf E }[V_iV_j].$$
We start with a bound for the summand ${\rm \bf E }[V_iV_j]$ when $j
\geq i+2$.
 First, notice that here
 \begin{equation}
{\rm \bf E }[V_iV_j]= {\rm \bf E }\Big[\Big(\sum_{h= \sigma^{i-1}+1}^{\sigma^i
}\frac {U_h}{h}\Big) \Big(\sum_{k= \sigma^{j-1}+1}^{\sigma^j }\frac
{U_k}{k}\Big)\Big]
 \end{equation}
 and
\begin{equation}\label{numerata}
h \leq \sigma^i \leq \sigma^{j-2} \leq \frac{k}{\sigma},\end{equation} 
hence, by assumption (i), $$g(h,k)\leqslant C;$$
thus the inequality in (v) can be simplified into
\begin{align*}
&|Cov(U_h, U_k)|\leqslant\frac{1}{h^\alpha} f\Big(h,\frac{k}{h}\Big)
\end{align*}  
(we incorporate the constant $C$ into $f$ for simplicity).
Hence, for $m$ sufficiently large in order that $\sigma^{j-1}+1 > m_0$, by assumption (ii) we have
\begin{align*}&{\rm \bf E }[V_iV_j]=\sum_{h= \sigma^{i-1}+1}^{\sigma^i
} \sum_{k= \sigma^{j-1}+1}^{\sigma^j }\frac
{1}{hk}{\rm \bf E }[U_hU_k]\leqslant
 \sum_{h= \sigma^{i-1}+1}^{\sigma^i
} \frac{1}{h^{\alpha+1}}\int_{ \sigma^{j-1}}^{\sigma^j } \frac
{1}{y}f\Big(h,\frac{v}{h}\Big)\, dv\\&\leqslant  \int_{ \sigma^{i-1}}^{\sigma^i
}du \,\frac{1}{u^{\alpha+1}} \int_{ \sigma^{j-1}}^{\sigma^j} \frac
{1}{v}f\Big(u,\frac{v}{u}\Big)\, dv.\end{align*}
By means of the change of variable $v =uz$ in the inner integral, the above becomes 
\begin{align*}
&\int_{ \sigma^{i-1}}^{\sigma^i
}dx \,\frac{1}{u^{\alpha+1}} \int_{ \frac{\sigma^{j-1}}{u}}^{\frac{\sigma^{j}}{u}} \frac
{1}{z}f(u,z)\, dz\leqslant\int_{ \sigma^{i-1}}^{\sigma^i
}dx \,\frac{1}{u^{\alpha+1}} \int_{ \frac{\sigma^{j-1}}{u}}^{\frac{\sigma^{j}}{u}} \phi(z)\, dz\\ &= \int_{ \sigma^{i-1}}^{\sigma^i
}du \,\frac{1}{u^{\alpha+1}}\Bigg\{F\Big(\frac{\sigma^{j}}{u}\Big)-F\Big(\frac{\sigma^{j-1}}{u}\Big)\Bigg\}.\end{align*}
Hence
\begin{align}\label{t}&\nonumber\sum_{m+1\leq i<j\leq
m+n \atop j
\geq i+2}{\rm \bf E }[V_iV_j]\leqslant \sum_{i=m+1} ^{m+n-2} \sum_{j=i+2} ^{m+n}\int_{ \sigma^{i-1}}^{\sigma^i
}du \,\frac{1}{u^{\alpha+1}}\Bigg\{F\Big(\frac{\sigma^{j}}{u}\Big)-F\Big(\frac{\sigma^{j-1}}{u}\Big)\Bigg\}\\&\nonumber=\sum_{i=m+1} ^{m+n-2}\int_{ \sigma^{i-1}}^{\sigma^i
}du \,\frac{1}{u^{\alpha+1}}\sum_{j=i+2} ^{m+n}\Bigg\{F\Big(\frac{\sigma^{j}}{u}\Big)-F\Big(\frac{\sigma^{j-1}}{u}\Big)\Bigg\}\\&\nonumber=\sum_{i=m+1} ^{m+n-2}\int_{ \sigma^{i-1}}^{\sigma^i
}du \,\frac{1}{u^{\alpha+1}}\Bigg\{F\Big(\frac{\sigma^{m+n}}{u}\Big)-F\Big(\frac{\sigma^{i+1}}{u}\Big)\Bigg\}\\&\nonumber=\sum_{i=m+1} ^{m+n-2}\int_{ \sigma^{i-1}}^{\sigma^i
}du \,\frac{1}{u^{\alpha+1}}F\Big(\frac{\sigma^{m+n}}{u}\Big)-\sum_{i=m+1} ^{m+n-2}\int_{ \sigma^{i-1}}^{\sigma^i
}du \,\frac{1}{u^{\alpha+1}}F\Big(\frac{\sigma^{i+1}}{u}\Big)\\&= \int_{ \sigma^{m}}^{\sigma^{m+n-2}
}du \,\frac{1}{u^{\alpha+1}}F\Big(\frac{\sigma^{m+n}}{u}\Big)-\sum_{i=m+1} ^{m+n-2}\int_{ \sigma^{i-1}}^{\sigma^i
}du \,\frac{1}{u^{\alpha+1}}F\Big(\frac{\sigma^{i+1}}{u}\Big).\end{align}

\noindent
By the change of variable $\frac{\sigma^{m+n}}{u}= v$ we get
\begin{align}\label{A}&\nonumber\int_{ \sigma^{m}}^{\sigma^{m+n-2}
}du \,\frac{1}{u^{\alpha+1}}F\Big(\frac{\sigma^{m+n}}{u}\Big)=\frac{1}{\sigma^{\alpha(m+n)}}\int_{ \sigma^{2}}^{\sigma^{n}
}F(v) v^{\alpha -1} \, dv\\&= \frac{1}{\sigma^{\alpha(m+n)}}\Big(\int_{\sigma}^{\sigma^{n}}F(v) v^{\alpha -1} \, dv-\int_{\sigma}^{\sigma^{2}}F(v) v^{\alpha -1} \, dv\Big).\end{align}
In a similar way, by the change of variable $\frac{\sigma^{i+1}}{u}= v$ we get
\begin{align} \label{B}&\int_{ \sigma^{i-1}}^{\sigma^{i}
}du \,\frac{1}{u^{\alpha+1}}F\Big(\frac{\sigma^{i+1}}{u}\Big)=\frac{1}{\sigma^{\alpha(i+1)}}\int_{ \sigma}^{\sigma^{2}
}F(v) v^{\alpha -1} \, dv\end{align}

\noindent
By inserting \eqref{A} and \eqref{B} into \eqref{t}, we get
\begin{align}
&\nonumber\sum_{m+1\leq i<j\leq
m+n \atop j
\geq i+2}{\rm \bf E }[V_iV_j]\\&\nonumber\leqslant\frac{1}{\sigma^{\alpha(m+n)}}\Big(\int_{\sigma}^{\sigma^{n}}F(u) u^{\alpha -1} \, du-\int_{\sigma}^{\sigma^{2}}F(u) u^{\alpha -1} \, du\Big)-\sum_{i=m+1} ^{m+n-2}\frac{1}{\sigma^{\alpha(i+1)}}\int_{ \sigma}^{\sigma^{2}
}F(u) u^{\alpha -1} \, du\\&\nonumber =\frac{1}{\sigma^{\alpha(m+n)}}\int_{\sigma}^{\sigma^{n}}F(u) u^{\alpha -1} \, du-\sum_{i=m+2} ^{m+n}\frac{1}{\sigma^{\alpha i}}\underbrace{\int_{ \sigma}^{\sigma^{2}
}F(u) u^{\alpha -1} \, du}_{= C}\\& =\frac{1}{\sigma^{\alpha(m+n)}}\int_{\sigma}^{\sigma^{n}}F(u) u^{\alpha -1} \, du-C\sum_{i=m+2} ^{m+n}\frac{1}{\sigma^{\alpha i}}.
\end{align}
And now, by \eqref{x}, \eqref{y} and \eqref{z},
 \begin{align}& \nonumber \label{fourth}\sum_{m+1\leq i<j\leq m+n}{\rm \bf E
}[V_iV_j]=\sum_{m+1\leq i<j\leq
m+n \atop j
\geq i+2} {\rm \bf E }[V_iV_j] +
\sum_{ i=m+1}^{m+n-1}{\rm \bf E }[V_iV_{i+1}]=
\\ \nonumber & \leqslant
\sum_{m+1\leq i<j\leq
m+n \atop j
\geq i+2}{\rm \bf E }[V_iV_j]
+
\sum_{ i=m+1}^{m+n-1}{\rm \bf E }[V_i^2]^{1/2}{\rm \bf E
}[V_{i+1}^2]^{1/2}
\\& \leqslant \frac{1}{\sigma^{\alpha(m+n)}}\int_{\sigma}^{\sigma^{n}}F(u) u^{\alpha -1} \, du-C\sum_{i=m+2} ^{m+n}\frac{1}{\sigma^{\alpha i}} +C \sum_{ i=m+1}^{m+n-1}1 \leqslant C\Big(n+\frac{1}{\sigma^{\alpha(m+n)}}\int_{\sigma}^{\sigma^{n}}F(u) u^{\alpha -1} \, du\Big).\end{align}
From \eqref{init}, \eqref{primaparte} and \eqref{fourth} we obtain
\begin{equation*}
{\rm \bf E}\Big[\big(\sum_{i=m+1}^{m+n }V_i\big)^2\Big]\leqslant 
C\Big(n+\frac{1}{\sigma^{\alpha(m+n)}}\int_{\sigma}^{\sigma^{n}}F(u) u^{\alpha -1} \, du\Big),\end{equation*}
i.e. the claim.

\hfill $\Box$

\noindent
Similar techniques prove the following more general result:
\begin{theorem}\label{secondorisultato}
Let $(U_n)_{n \geqslant 1}$ be a sequence of centered random variables. 
Let $N$ be an integer and assume that there exist a number $\sigma >1$ and for each $j=1, 2, \dots, N$ numbers $\alpha_j \geqslant 0$  a non--negative function  $f_j(u,z)$ defined on the set $\{u \geqslant 1, z \geqslant \sigma \}$, a non--negative double--indexed sequence $g_j$ defined on the set $\{(m,n)\in \mathbb{N}^2: \sigma m\leqslant n\}$ such that
$$
\sup_{n \geqslant \sigma m} g_j(m,n) = C < + \infty; \leqno \it(i)$$
(ii) uniformly in $u>0$ the functions $u \mapsto f_j\big(u, \frac{v}{u}\big)$ are ultimately non--increasing  (i.e.  there exists $m_0$ such that $v \mapsto f_j\big(u, \frac{v}{u}\big)$  are non--increasing on $(m_0, + \infty)$, for each $j=1, \dots, N$ and for every $u>0$);

\noindent
(iii) for each $j= 1, 2, \dots, N$ the functions $$z \mapsto \phi_j(z) = \sup_{u \geqslant 1} \frac{f_j(u,z)}{z},\qquad u \mapsto F_j(u)=\int_\sigma^ u\phi_j(z) \, dz $$
are defined on $[\sigma, + \infty)$;
$$
|Cov(U_m, U_n)|\leqslant C\begin{cases}
m & \hbox{\rm for } m=n\\1& \hbox{\rm for } m<n \leqslant \sigma m;
\end{cases}\leqno \it(iv)
$$
(v) there exists $m_1$ such that, for $n> m \geqslant m_1$
$$
|Cov(U_m, U_n)|\leqslant \sum_{j=1}^N g_j(m,n)\frac{1}{m^{\alpha_j}} f_j\Big(m,\frac{n}{m}\Big).
$$
Denote
$$V_n = \sum_{k=\sigma^{n-1}+1}^{\sigma^{n}} \frac{U_k}{k}.$$
Then, for every $n$ and every sufficiently large $m$
\begin{equation*}
E\big[(\sum_{i=m+1}^{m+n}V_i)^2\big]\leqslant C\left(n + \sum_{j=1}^N\frac{1}{\sigma^{\alpha_j (m+n)}}\int_{\sigma}^{\sigma^n} F_j(u) u^{\alpha_j -1}\, du\right).
\end{equation*}

\end{theorem}

\begin{corollary} \label{primocorollario}  In the setting of Theorem \ref{primoriusltato}, assume in addition that $\alpha =0$ and there exists $\beta >1$ such that $F(x) \leqslant C (\log x)^{\beta -1}$ for every $x > \sigma$. Then, for every sufficiently large $m$, \begin{equation}
{\rm \bf E}\Big[\big(\sum_{i=m+1}^{m+n }V_i\big)^2\Big]\leqslant 
C\Big((m+n)^\beta  - m^\beta\Big).\end{equation}\end{corollary}
\bigskip
{\it Proof.} Putting $\alpha =0$ in the claim of Theorem \ref{primoriusltato} we obtain
\begin{align*}
&{\rm \bf E}\Big[\big(\sum_{i=m+1}^{m+n }V_i\big)^2\Big]\leqslant 
C\Big(n+\int_{\sigma}^{\sigma^{n}}\frac{F(u)}{u} \, du\Big)\leqslant C\Big(n+\int_{\sigma}^{\sigma^{n}}\frac{(\log u)^{\beta -1}}{u} \, du\Big)=C\Big(n+\Big[(\log u)^{\beta}\Big]_\sigma^{\sigma^{n}}\Big)\\&\leqslant C\Big(n+n^{\beta}\Big)\leqslant C n^\beta.
\end{align*}
On the other hand, the function $z\mapsto \Big\{(z+n)^\beta-z^\beta\Big\}$ being increasing (its derivative is $\beta(z+n)^{\beta -1}-\beta z^{\beta -1}\geqslant 0 $), we have $$n^\beta \leqslant(1+n)^\beta-1 \leqslant (m+n)^\beta-m^\beta.$$
\hfill $\Box$
\begin{corollary} \label{secondocorollario} In the setting of Theorem \ref{primoriusltato}, assume in addition that $\alpha >0$ and there exists $\beta >1$ such that $F(x) \leqslant C (\log x)^{\beta }$ for every $x > \sigma$. Then, for every sufficiently large $m$, \begin{equation}
{\rm \bf E}\Big[\big(\sum_{i=m+1}^{m+n }V_i\big)^2\Big]\leqslant 
C\Big((m+n)^\beta  - m^\beta\Big).\end{equation}\end{corollary}
\bigskip
{\it Proof.} In this case Theorem \ref{primoriusltato} gives\begin{align*}&
{\rm \bf E}\Big[\big(\sum_{i=m+1}^{m+n }V_i\big)^2\Big]\leqslant 
C\Big(n+\frac{1}{\sigma^{\alpha(m+n)}}\int_{\sigma}^{\sigma^{n}}F(u) u^{\alpha -1} \, du\Big)\leqslant 
C\Big(n+\frac{1}{\sigma^{\alpha n}}\int_{\sigma}^{\sigma^{n}}(\log u)^{\beta } u^{\alpha -1} \, du\Big)\\& \leqslant C\Big(n+\frac{1}{\sigma^{\alpha n}}\big(\log (\sigma^n)\big)^{\beta } \int_{\sigma}^{\sigma^{n}}u^{\alpha -1} \, du\Big)=C\Big(n+\frac{\big(\log (\sigma^n)\big)^{\beta }(\sigma^{\alpha n}- \sigma^\alpha)}{\alpha\sigma^{\alpha n}}\Big)\leqslant C\Big(n+n^{\beta}\Big)\leqslant C n^\beta.\end{align*}
The rest of the proof is identical to Corollary \ref{primocorollario}.

\hfill $\Box$

\begin{corollary} \label{terzocorollario} In the setting of Theorem \ref{secondorisultato}, assume that there exists $\beta >1$ such that $\sum_{j=1}^NF_j(x) \leqslant C(\log x)^{\beta }$ for every $x > \sigma$. Then

\begin{itemize}
\item[(i)]  for every sufficiently large $m$ and for every $n$, \begin{equation*}
{\rm \bf E}\Big[\big(\sum_{i=m+1}^{m+n }V_i\big)^2\Big]\leqslant 
C\Big((m+n)^\beta  - m^\beta\Big).\end{equation*}
\item[(ii)]for every $\delta > 0$, $$\sum_{ i=1}^n V_i =
O( n^{\beta/2}(\log n)^{2+\delta} ), \quad P- a.s.$$

\end{itemize}
\end{corollary}

\noindent
{\it Proof.} Point (i) follows from Corollaries \ref{primocorollario} and \ref{secondocorollario}. Point (ii) is a consequence of the well known Gaal--Koksma Strong Law of Large Numbers (see
\cite{PS}, p. 134); here is the precise statement: \begin{theorem}
\label{GK}Let $(V_n)_{n\geqslant 1}$ be a sequence of centered random variables
with finite variance. Suppose that there exists a constant $\beta >
0$ such that, for all integers $m \geq 0$, $n \geqslant 1 $, \begin
{equation}\label{bound}
E\Big[\big(\sum_{i=m+1}^{m+n}V_i\big)^2\Big]\leq C\big((m+n)^\beta-
m^\beta\big),\end {equation} for a suitable constant $C$ independent
of $m$ and $n$. Then, for every $\delta > 0$, $$\sum_{ i=1}^n V_i =
O( n^{\beta/2}(\log n)^{2+\delta} ), \quad P- a.s.$$
\end{theorem}

\begin {remark} \label{ultimate} \rm It is not difficult to see that Theorem \ref{GK}
 is in force even if the bound
\eqref{bound} holds only for all integers $m \geq h_0$, $n> 0$,
where $h_0$ is an integer strictly greater than 0. A rigorous proof of this statement can be found in the appendix of \cite{GAS}. From now on, this slight generalization will be tacitly used.
\end {remark}

\noindent
\begin{theorem}({\bf General ASLT}) \label{principale} 
Let $(Y_n)_{n\geqslant 1}$ be a sequence of non--negative (resp. non--positive) random variables with
$$\lim_{n \to \infty}E[Y_n] = \ell >0\qquad (\hbox{resp. }\ell <0)$$
and such that the sequence $(U_n)_{n\geqslant 1}$ defined by $U_n = Y_n-E[Y_n]$ verifies the assumptions of Theorem \ref{secondorisultato}. Assume that there exists $\beta >1$ such that $\sum_{j=1}^NF_j(x) \leqslant C(\log x)^{\beta }$ for every $x > \sigma$. Then
$$\lim_{n \to \infty}\frac{1}{\log n}\sum_{k=1}^n \frac{Y_k}{k}=\ell, \qquad a.s.$$

\end{theorem}  
{\it Proof.} By point (ii) of Corollary \ref{terzocorollario}, for every $\delta >0$ we have
\begin{equation}\label{limite}
\frac{{\sum_{ i=1}^n V_i}}{n} =
\frac {O( n^{\beta/2}(\log n)^{2+\delta} )}{n} \mathop {\longrightarrow}_{n \to \infty} 0.
\end{equation}
Since
\begin{equation*}
\sum_{ i=1}^n V_i=\sum_{i=1}^n \sum_{k=\sigma^{i-1}+1}^{\sigma^{i}} \frac{U_k}{k}=\sum_{k=2}^{\sigma^{n}} \frac{U_k}{k}= \sum_{k=2}^{\sigma^{n}} \frac{Y_k}{k}-\sum_{k=2}^{\sigma^{n}} \frac{E[Y_k]}{k}
\end{equation*}
and
\begin{equation*}
\frac{1}{n \log \sigma}\sum_{k=2}^{\sigma^{n}} \frac{E[Y_k]}{k}\mathop {\longrightarrow}_{n \to \infty} \ell,
\end{equation*}
the relation \eqref{limite} is equivalent to
\begin{equation*}
 \frac{1}{n \log \sigma}\sum_{k=2}^{\sigma^{n}} \frac{Y_k}{k}\mathop {\longrightarrow}_{n \to \infty} \ell;
\end{equation*}
By the same argument as in \cite{GAS}, pp. 789--790,  this in turn implies that
\begin{equation*}
\frac{1}{\log n}\sum_{k=1}^{ n} \frac{Y_k}{k}\mathop {\longrightarrow}_{n \to \infty} \ell,
\end{equation*}
i.e. the claim.

\hfill $\Box$

\section{The Almost Sure Local Limit Theorem}
Let $x>0$ be given; let $\boldsymbol \kappa =(\kappa_n )_{n\geqslant 1}$ be a strictly increasing sequence of integers with $\kappa_n \sim x n$, fixed throughout the sequel. Let $(Y_n)_{n\geqslant 1}$ be the sequence defined in \eqref{definizione}; the main result of this section and of the paper (Theorem \ref{ASLLT}) is an ASLLT for the sequence $(Y_n)_{n\geqslant 1}$. Before proving it, we need some Lemmas.  
  
  \smallskip
\noindent
For every $\epsilon \in (0, \frac{1}{2x}) $ we set
 \begin{equation}\label{sigma}
 \sigma=\sigma_{\epsilon}= \frac{1+x(1-\epsilon)}{x(1+\epsilon)}=1+ \frac{2}{1+\epsilon}\Big(\frac{1}{2x}- \epsilon \Big)>1.
 \end{equation}
 \begin{lemma}\label{lemmino} Let $\epsilon\in (0, \frac{1}{2x}) $ be fixed. Then there exists $m_0 =m_0(\epsilon )$ such that, for $\sigma m> n >m > m_0$,
 $$P\big(T_m^n =  \kappa_n-\kappa_m\big)=0$$ \end{lemma}
 \bigskip
 \noindent
 {\it Proof.} Let 
$$A = \bigcap_{k=m+1}^n \{Z_k =0\}.$$
Then
$$P\big(T_m^n =  \kappa_n-\kappa_m\big)= P\big(\{T_m^n =  \kappa_n-\kappa_m\}\cap A\big)
+P\big(\{T_m^n =  \kappa_n-\kappa_m\}\cap A^c\big).$$

\begin{itemize}
\item[(i)] Let $m_0$ be such that, for every $m> m_0$,
$$xm(1-\epsilon) < \kappa_m <xm(1+\epsilon). $$Then, for $\sigma m> n >m > m_0$,
\begin{align}\label{HH}&\nonumber\kappa_n-\kappa_m< xn(1+\epsilon)- xm (1-\epsilon)=m\Big\{x \frac{n}{m} (1+\epsilon) - x(1-\epsilon)\Big\}\\&\leqslant m\Big\{x \sigma (1+\epsilon) - x(1-\epsilon)\Big\}=m.
\end{align}
Hence
\begin{align*}
&\{T_m^n = \kappa_n-\kappa_m\}\cap A^c =\{T_m^n = \kappa_n-\kappa_m\}\cap
\Big( \bigcup_{k=m+1}^n \{Z_k =1\} \Big)\cr &=\bigcup_{k=m+1}^n
\{T_m^n = \kappa_n-\kappa_m, Z_k =1\} \subseteq \bigcup_{k=m+1}^n  \{T_m^n =
\kappa_n-\kappa_m, T_m^n\geqslant m+1\}\\ & = \{T_m^n = \kappa_n-\kappa_m, T_m^n\geqslant
m+1\}= \emptyset,
\end{align*}
by \eqref{HH}.
\item[(ii)]$$A \subseteq  \{ T_m^n=0\},$$ hence 
\begin{align*}
&\{T_m^n =\kappa_n-\kappa_m \}\cap A \subseteq \{T_m^n = \kappa_n-\kappa_m, T_m^n=0\}= \emptyset,\end{align*}

since $ \kappa_n-\kappa_m> 0$.
\end{itemize}

\hfill $\Box$

  \begin{lemma}\label{lemmaausiliario} Let $\epsilon\in (0, \frac{1}{2x}) $ be fixed. Then there exists $m_0 =m_0(\epsilon )$ such that, for $n \geqslant m > m_0$,
$$|Cov(Y_m, Y_n)|\leqslant C\begin{cases}
   m & {\rm for}\,\,  m=n\\1 & {\rm for}\,\,  m<n \leqslant \sigma m,
  \end{cases} $$ 
  where $C$ is a positive constant.
\end{lemma}

  {\it Proof.} 
  
  \begin{itemize}
  \item[(a)]For $m=n$:\begin{align*}
&Cov(Y_m, Y_m)= m^2\big\{ P(T_m = \kappa_m )-P^2(T_m =\kappa_m  )\big\}\\&=\big\{mP(T_m =\kappa_m )\big\}\big\{m-P(T_n = \kappa_m )\big\}\leqslant  Cm,
\end{align*}
 by the Local Theorem (Corollary \ref{local}).
  
  \item[(b)] For $ m< n\leqslant\sigma m$: let $m_0$ be as in Lemma \ref{lemmino}. Then, for $\sigma m>n>m> m_0$,
  \begin{align*}
  &|Cov(Y_m, Y_n)|= mn \big| P(T_m = \kappa_m,T_n = \kappa_n )-P(T_m = \kappa_m)P(T_n = \kappa_n )\big|\\&=\big\{mP(T_m = \kappa_m)\big\}
\big|n P\big(T_m^n= \kappa_n-\kappa_m\big) -n P(T_n = \kappa_n)\big|\\&=\big\{mP(T_m = \kappa_m)\big\}
\big\{n P(T_n =\kappa_ n)\big\}\leqslant C,\end{align*}
by the Local Theorem again and observing that
$ P\big(T_m^n= \kappa_n-\kappa_m \big)= 0$, by Lemma \ref{lemmino}.

  \end{itemize}

\hfill $\Box$
 
 \begin{remark}\rm  Notice that 
 \begin{itemize}
 \item[(i)] $\kappa_n=[xn]$ is strictly increasing if $x \geq 1$;
 \item[(ii)] If $\kappa_n=xn$ we can take $\epsilon =0$ and $m_0=1$.
 \end{itemize}

 \end{remark}
 
\begin{lemma} \label{lemma}  In the setting of Theorem \ref{primoriusltato}, assume that $f$ has the form
$$f(u,z) = \psi(uz)$$
where $t \mapsto \psi(t)$ is a continuous ultimately non--increasing function, i.e. there exists $t_0$ such that $t \mapsto \psi(t)$ is non--increasing for $t \geqslant t_0$. Then
 \begin{equation}
F(u) \leqslant C\begin{cases}1 & \hbox{ for u}\leqslant t_0\\
1 + \int_{t_0}^u \frac{\psi(z)}{z}\, dz& \hbox{ for u}> t_0.
\end{cases} \end{equation}\end{lemma}
\bigskip
{\it Proof.} It is easy to see that
\begin{equation*}
\sup_{x \geqslant 1}f(x,z) =\sup_{x \geqslant 1}\psi(xz)=
\sup_{u \geqslant z}\psi(u)\begin{cases}\leqslant \displaystyle \max_{u\in [1,t_0]}\psi(u)=:M& \hbox{ for } z \leqslant t_0\\
=\psi(z)& \hbox{for } z > t_0.
\end{cases}
\end{equation*}
Hence $$\phi(z) \leqslant\begin{cases} \frac{M}{z}& \hbox{for } z \leqslant t_0\\
\frac{\psi(z)}{z}& \hbox{for } z > t_0.\end{cases}$$
and

\begin{equation*}
F(u) =\int_\sigma^u \phi(z) \, dz = \underbrace{\int_\sigma^{t_0}\phi(z) \, dz}_{=C}+ \int_{t_0}^u\phi(z)\, dz \leqslant C\begin{cases}1 & \hbox{ for }u\leqslant t_0\\
1 + \int_{t_0}^u \frac{\psi(z)}{z}\, dz& \hbox{ for }u> t_0.
\end{cases} 
\end{equation*}
\hfill $\Box$
\begin{remark}
\rm Of course, the preceding lemma has an obvious generalization in the setting of Theorem \ref{secondorisultato}.
\end{remark}

\bigskip
\noindent
We are ready to give the

 \bigskip
 \noindent
 {\it Proof of Theorem \ref{ASLLT}.} Though with tedious and cumbersome calculations, it is easy to see that the correlation inequality of Theorem \ref{covarianzabase} takes a slightly more tractable form for sufficiently large $m$: precisely (we neglect the multiplicative constant $C$ for easy writing):
  \begin{align}\label{formaalternativa}
 &\nonumber |Cov(Y_m, Y_n)|\leqslant  \frac{n}{\kappa_n - \kappa_m}\Bigg\{ \frac{1+\log\frac{n}{m}}{\sqrt{n-m}}+ \Big[\exp \Big(C \Big\{\frac{\log^3\frac{n}{m}}{(n-m)^2}+ \frac{m \log 
  ^2\frac{n}{m}}{n-m}\Big\}\Big)-1 + \frac{1}{\log \frac{n}{m}}\Big]\\&+ \nonumber\frac{\big|(xn-\kappa_n)-(xm-\kappa_m)\big|}{n-m} + \frac{m}{n-m}\Bigg\}+  \frac{m}{n-m}\\&+
 \frac{\log n}{\sqrt{n}} +\exp \Big(C \Big\{\frac{\log^3 n}{n^2}+ \frac{ \log 
  ^2n}{n}\Big\}\Big)-1 + \frac{1}{\log n}+\frac{\big|(xn-\kappa_n)-(2x-\kappa_2)\big|}{n-2}+ \frac{1}{n}.\end{align}

\noindent  
(In fact (look at the formula in the statement of Theorem \ref{covarianzabase})
\begin{align*}
\frac{n}{n-m}\chi_{m,n}^{\boldsymbol \kappa, x}=\underbrace{\frac{n}{\kappa_n-\kappa_m}\Big(\frac{\log \frac{n}{m}}{\sqrt{n-m}}+g_{m,n}\Big)}_{(a)}+\underbrace{x\Big|\frac{n-m}{\kappa_n-\kappa_m}- \frac{1}{x}\Big|\cdot\frac{n}{n-m}}_{(b)}+\underbrace{\frac{n}{n-m}\cdot \frac{m+1}{\kappa_n-\kappa_m}}_{(c)},
\end{align*}
  and
  \begin{itemize}
  \item[(a)]
  \begin{align*}
  & \frac{n}{\kappa_n-\kappa_m}\Big(\frac{1+\log \frac{n}{m}}{\sqrt{n-m}}+g_{m,n}\Big)\\&= \frac{n}{\kappa_n-\kappa_m}\Bigg[\frac{1+\log \frac{n}{m}}{\sqrt{n-m}}+\exp\Big(C \Big\{\frac{\log \frac{n}{m}}{(n-m)^2}+\frac{m+2}{n-m}\Big\}\log^2 \frac{n}{m}\Big)-1 + \frac{1}{\log \frac{n}{m}}\Bigg]\\&\leqslant \frac{n}{\kappa_n-\kappa_m}\Bigg[\frac{1+\log \frac{n}{m}}{\sqrt{n-m}}+\exp\Big(C \Big\{\frac{\log^3 \frac{n}{m}}{(n-m)^2}+\frac{m\cdot\log^2 \frac{n}{m}}{n-m}\Big\}\Big)-1 + \frac{1}{\log \frac{n}{m}}\Bigg] ;\end{align*}
  \item[(b)]
   \begin{align*}x\Big|\frac{n-m}{\kappa_n-\kappa_m}- \frac{1}{x}\Big|\cdot\frac{n}{n-m}=\frac{n}{\kappa_n-\kappa_m}\cdot\Big|\frac{(nx-\kappa_n)-(mx-\kappa_m)}{n-m}\Big|;
  \end{align*}
  \item[(c)]\begin{align*}\frac{n}{n-m}\cdot \frac{m+1}{\kappa_n-\kappa_m}\leqslant\frac{n}{\kappa_n-\kappa_m}\cdot\frac{m}{n-m}.\end{align*}
  \end{itemize}
  Further, by (a), (b) and (c) above
   \begin{align*}& \chi_{2,n}^{\boldsymbol \kappa, x}\leqslant \frac{n}{\kappa_n - \kappa_2}\Bigg[\frac{1+\log\frac{n}{2}}{\sqrt{n-2}}+ \exp\Big(C\Big\{\frac{\log^3 \frac{n}{2}}{(n-2)^2}+\frac{2\cdot\log^2 \frac{n}{2}}{n-2}\Big\}\Big)-1+\frac{1}{\log \frac{n}{2}}+\\&+\frac{n}{n-2}\cdot \frac{(nx-\kappa_n)-(2x-\kappa_2)}{n-2}+\frac{1}{n-2}\Bigg]\\&\leqslant\frac{\log n}{\sqrt n}+\exp \Big(C \Big\{\frac{\log^3 n}{n^2}+ \frac{ \log 
  ^2n}{n}\Big\}\Big)-1 + \frac{1}{\log n}+\frac{\big|(xn-\kappa_n)-(2x-\kappa_2)\big|}{n-2}+ \frac{1}{n} ,\end{align*}
   for sufficiently large $n$; recall that we are neglecting multiplicative constants).
   
   \bigskip
   \noindent
   The statement of the Theorem is a consequence of the general ASLT (Theorem \ref{principale}): we check  assumption (v) of Theorem \ref{secondorisultato} for each summand in the basic correlation inequality (in the form \eqref{formaalternativa}) and use Corollary \ref{primocorollario} or Corollary \ref{secondocorollario}, as needed (it is easy to see that assumptions (i)--(iv) of Theorem \ref{secondorisultato} are in force for each summand in the basic correlation inequality, hence we omit the details). Precisely (with the notations of Theorem \ref{secondorisultato} and with $\sigma$ defined in \eqref{sigma}):

\begin{itemize}

\item[(1)] First summand: $$\frac{n}{\kappa_n - \kappa_m}\cdot\frac{1+\log\frac{n}{m}}{\sqrt {n-m}}.$$ 
Fix $\delta \in \Big(0, \frac{\sigma -1}{\sigma +1}\Big)$, and let $m_1$ be such that $$1-\delta< \frac{\kappa_n}{n}< 1+\delta, \qquad n > m_1.$$
 Then, for $n> \sigma m$ and $m> m_1$,\begin{equation}
 \frac{\kappa_n - \kappa_m}{n}= \frac{\kappa_n}{n}- \frac{\kappa_m}{m}\cdot \frac{m}{n}\geqslant (1-\delta)- \frac{1+\delta}{\sigma}>0,
 \end{equation}
 hence
 \begin{equation}
 \sup_{n> \sigma m}\frac{n}{\kappa_n - \kappa_m}\leqslant \frac{1}{ (1-\delta)- \frac{1+\delta}{\sigma}}.
 \end{equation}
Moreover we have
$\frac{1+\log\frac{y}{x}}{\sqrt {y-x}}=\frac{1}{\sqrt x}\cdot\frac{1+\log\frac{y}{x}}{\sqrt {\frac{y}{x}-1}}$, hence
$$g_1(m,n) =\frac{n}{\kappa_n - \kappa_m}, \qquad \alpha_1 = \frac{1}{2}, \quad f_1(u,z) =\frac{1+\log z}{\sqrt {z-1}}, \qquad  \phi_1(z)=\frac{1+\log z}{z\sqrt {z-1}};$$

last
$$F_1(u)= \int_\sigma^u\phi_1(z)\, dz = \int_\sigma^u\frac{1+\log z}{z\sqrt {z-1}}\, dz \leqslant C \leqslant C( \log u)^\beta, \quad \forall \beta >1.$$

\item[(2)]Second summand: $$\frac{n}{\kappa_n - \kappa_m}\cdot\Big\{\exp\Bigg(C\Big\{\frac{\log^3\frac{n}{m}}{(n-m)^2}+ \frac{m\log^2\frac{n}{m}}{n-m}\Big\}\Bigg)-1\Big\}.$$
We have again $g_2(m,n) =\frac{n}{\kappa_n - \kappa_m}$; moreover$$\exp\Bigg(C\Big\{\frac{\log^3\frac{y}{x}}{(y-x)^2}+ \frac{x\log^2\frac{y}{x}}{y-x}\Big\}\Bigg)-1=\exp\Bigg(C\Big\{\frac{\log^3\frac{y}{x}}{x^2(\frac{y}{x}-1)^2}+ \frac{x\log^2\frac{y}{x}}{x(\frac{y}{x}-1)}\Big\}\Bigg)-1,$$
so that 
$$f_2(u,z) = \exp\Bigg(C\Big\{\frac{\log^3 z}{u^2(z-1)^2}+ \frac{x\log^2 z}{u(z-1)}\Big\}\Bigg)-1,$$
and $\alpha_2 =0$; further $$\phi_2(z) =\sup_{u \geqslant 1}\frac{f_2(u,z)}{z}=\frac{1}{z}\Bigg\{\exp\Bigg(C\Big\{\frac{\log^3 z}{(z-1)^2}+ \frac{\log^2z }{z-1}\Big\}\Bigg)-1\Bigg\}.$$
Put
$$M = \sup_{z \geqslant \sigma}\frac{\exp\Bigg(C\Big\{\frac{\log^3 z}{(z-1)^2}+ \frac{\log^2 z}{z-1}\Big\}\Bigg)-1}{\Big\{\frac{\log^3 z}{(z-1)^2}+ \frac{\log^2 z}{z-1}\Big\}}$$
Then $$\phi_2(z)\leqslant M \Big\{\frac{\log^3 z}{(z-1)^2}+ \frac{\log^2 z}{z-1}\Big\}\frac{1}{z}$$
and $$F_2(u)= \int_\sigma^u\phi(z)\, dz \leqslant M\int_\sigma^u \Big\{\frac{\log^3 z}{(z-1)^2}+ \frac{\log^2 z}{z-1}\Big\}\frac{1}{z}\, dz \leqslant C \leqslant C( \log u)^\beta, \quad \forall \beta >1.$$
\item[(3)]Third summand: $$\frac{n}{\kappa_n-\kappa_m}\cdot\frac{1}{\log\frac{n}{m}}$$
we have $g_3(m,n) =\frac{n}{\kappa_n - \kappa_m}$, $\alpha_3=0$ and 
$$ f_3(u,z)= \frac{1}{\log z}; \qquad \phi_3(z)= \frac{1}{z\log z},$$
hence
$$F_3(u)= \int_\sigma^u\phi_3(z)\, dz =\int_\sigma^u\frac{1}{z\log z}\,dz = \Big[\log \log z\Big]_\sigma^u\leqslant \log \log u \leqslant  (\log u)^\beta, \quad \forall \beta >1.$$
\item[(4)] Fourth summand: $$\frac{n}{\kappa_n-\kappa_m}\cdot\frac{\big|(xn-\kappa_n)-(xm-\kappa_m)\big|}{n-m}$$
Once more, $g_4(m,n) =\frac{n}{\kappa_n - \kappa_m}$, $\alpha_4=0$. Let $\delta > 0$ be fixed and $m_0$ such that $$|\kappa_n - nx|< \delta x n, \qquad n>m_0.$$
Then, for $n>m> m_0,$
$$\frac{\big|(xn-\kappa_n)-(xm-\kappa_m)\big|}{n-m}<\delta x\frac{n+m}{n-m}=\delta x\frac{\frac{n}{m}+1}{\frac{n}{m}-1}$$
and $$f_4\big(u,z\big)=\delta x\frac{z+1}{z-1}; \qquad \phi_4(z)= \delta x\frac{z+1}{z(z-1)}< \frac{C}{z}.$$
Hence
$$F_4(u)= \int_\sigma^u\phi(z)\, dz <C\int_\sigma^u\frac{1}{z}\, dz< C \log u \leqslant C \log^\beta u, \qquad \forall \, \beta >1.$$

\item[(5)] Fifth summand: 
$$\frac{n}{\kappa_n-\kappa_m}\cdot\frac{m}{n-m}= \frac{n}{\kappa_n-\kappa_m}\cdot \frac{1}{\frac{n}{m}-1}.$$
Once more, $g(m,n) =\frac{n}{\kappa_n - \kappa_m}$, $\alpha=0$ and $$f(u,z) = \frac{1}{z-1}; \qquad \phi(z) = \frac{1}{z(z-1)},$$
and$$F(u)= \int_\sigma^u\phi(z)\, dz =\int_\sigma^u\frac{1}{z(z-1)}\,dz \leqslant C \leqslant  (\log u)^\beta, \quad \forall \beta >1.$$
\item[(6)] Sixth summand:

$$\frac{m}{n-m}=  \frac{1}{\frac{n}{m}-1}.$$
Here $g_6(m,n) =1$, $\alpha_6=0$ and $$f_6(u,z) = \frac{1}{z-1}.$$
The argument is identical to the previous one.

\item[(7)] Seventh summand:
\begin{align*}
&\frac{\log n}{\sqrt{n}}+ \exp \Big(C \Big\{\frac{\log^3 n}{n^2}+ \frac{ \log 
  ^2n}{n}\Big\}\Big)-1 + \frac{1}{\log n}\\&\leqslant C\Big[\frac{\log n}{\sqrt{n}}+ \exp \Big(C \Big\{\frac{ \log 
  ^2n}{n}\Big\}\Big)-1 + \frac{1}{\log n}\Big].
\end{align*}
Here $g_7\equiv 1$, $\alpha_7=0$ and 
$$f_7(u,z) = \frac{\log uz}{\sqrt{uz}}+ \exp \Big(C \Big\{\frac{ \log 
  ^2uz}{uz}\Big\}\Big)-1 + \frac{1}{\log uz}=\psi(uz),$$
 with
 $$\psi_7(t)=\frac{\log t}{\sqrt{t}}+ \exp \Big(C \Big\{\frac{ \log 
  ^2t}{t}\Big\}\Big)-1 + \frac{1}{\log t}.$$
 
 We can apply Lemma \ref{lemma}, and we find 
$$F_7(u)= C + \int_{t_0}^u\Bigg(\frac{\log t}{t\sqrt{t}}+ \frac{1}{t}\Big\{\exp \Big(C \Big\{\frac{ \log 
  ^2t}{t}\Big\}\Big)-1\Big\} + \frac{1}{t\log t} \Bigg)\, dt\leqslant C \log\log u \leqslant (\log u)^\beta, \qquad \forall \beta >1,$$
  for some suitable $t_0> \sigma$.

\item[(8)] Eighth summand:
$$\frac{\big|(xn-\kappa_n)-(2x-\kappa_2)\big|}{n-2}\leqslant C.$$
Here
$g_8\equiv 1$, $\alpha_8=0$ and 
$$f_8(u,z) = C=\psi_8(uz),$$
 with $\psi_8(t)=C$.
 We can apply Lemma \ref{lemma}, and we find 
$$F_8(u)= C + \int_{t_0}^u\frac{1}{t} \, dt\leqslant C\log u  \leqslant (\log u)^\beta, \qquad \forall \beta >1,$$for some suitable $t_0> \sigma.$
\item[(9)] Ninth summand:$$\frac{1}{n}$$
The argument is the same as in (7) and (8).

\end{itemize}
 \hfill $\Box$

\section {Explicit form of the cumulants of the Bernoulli distribution
}

In this section we prove the explicit formula announced in  Remark \ref{cumulants}.
For every integer $n$ and every integer $k$ with $0 \leqslant k \leqslant n$ put
$$a_{k,n}= \sum_{j=0}^k (-1)^{j+1} {k\choose j}j^n. $$
\begin{remark}\label{1}
(i) Notice that $a_{1,n}=1$ for every $n$.

(ii) The Stirling number of second kind $S(n,k)$ has the explicit expression
$$S(n,k)= \frac{1}{k! } \sum_{j=0}^k (-1)^{k-j} {k\choose j}j^n.$$
Hence
\begin{align*}
& a_{k,n}= \sum_{j=0}^k (-1)^{j+1} {k\choose j}j^n= (-1)^{k+1} \sum_{j=0}^k (-1)^{j-k} {k\choose j}j^n= (-1)^{k+1} \sum_{j=0}^k (-1)^{k-j} {k\choose j}j^n\\ &=(-1)^{k+1}k!S(n,k).
\end{align*}
(ii) We also recall that $S(n,n)=1$, which  implies that $a_{n,n}=(-1)^{n+1}n!$ by the above relation.
\end{remark}
Let $\mathcal{B}(1,x)$ be the Bernoullian law with parameter $x \in (0,1)$. Denote by $c_n(x)$ the $n-$th cumulant of $\mathcal{B}(1,x)$, i.e. the $n-$th coefficient in the development of the logarithm of its  characteristic function $\phi(t)$:
$$\log \phi(t) = \log \Big(1+ x (e^{it}-1)\Big)= \sum_{n=1}^\infty
c_n(x) \frac{(it)^n}{n!}.$$
\begin{remark}\label{2a} (i) It is easily seen that $c_1(x) =x.$

(ii) It is well known (see \cite{K} ex. 6 p. 312 for instance) that the sequence of functions $\big(c_n(x)\big)_n$  verifies the recurrence relation
\begin{equation}\label{ricorrenza}
c_{n+1}(x) = x(1-x)c^\prime_{n}(x)
\end{equation}

\end{remark}

\begin{proposition}
For every $n \geqslant 2$ we have
\begin{equation}\label{a}
c_n(x) =x(1-x)\Big\{\sum_{k=1}^{n-1}a_{k,n-1}x^{k-1}\Big\}.
\end{equation}

\end{proposition}

{\emph Proof.} By \eqref{ricorrenza}, we must prove that, for every $n \geqslant 1$, 
\begin{equation}\label{dadimostrare1}
c^\prime_{n}(x)=\sum_{k=1}^{n}a_{k,n}x^{k-1}.
\end{equation}
The proof is by induction. 

 \noindent
For $n=1$ the statement follows from Remarks \ref{1} (i) and \ref{2a} (i). 

 \noindent
Assume that \eqref{dadimostrare1} holds for the integer $n-1$; hence, by \eqref{ricorrenza}, we have
\begin{equation*}
c_{n}(x)=x(1-x)\Big\{\sum_{k=1}^{n-1}a_{k,n-1}x^{k-1}\Big\}
\end{equation*}
and differentiating we get
 \begin{align*}
&c^\prime_{n}(x)=(1-2x)\Big\{\sum_{k=1}^{n-1}a_{k,n-1}x^{k-1}\Big\}+ (x-x^2)\Big\{\sum_{k=2}^{n-1}(k-1)a_{k,n-1}x^{k-2}\Big\}\\ &= 
\sum_{k=1}^{n-1}a_{k,n-1}x^{k-1}- \sum_{k=1}^{n-1}2a_{k,n-1}x^{k}+ \sum_{k=2}^{n-1}(k-1)a_{k,n-1}x^{k-1}-\sum_{k=2}^{n-1}(k-1)a_{k,n-1}x^{k}\\&=\sum_{k=1}^{n-1}a_{k,n-1}x^{k-1}-
\sum_{k=2}^{n}2a_{k-1,n-1}x^{k-1}+\sum_{k=2}^{n-1}(k-1)a_{k,n-1}x^{k-1}-\sum_{k=3}^{n}(k-2)a_{k-1,n-1}x^{k-1}\\& =a_{1,n-1}+\big(a_{2,n-1}-2a_{1,n-1}+a_{2,n-1}\big)x+\\&+ \sum_{k=3}^{n-1}x^{k-1}\Big(a_{k,n-1}-2a_{k-1,n-1}+(k-1)a_{k,n-1}
-(k-2)a_{k-1,n-1}\Big)+\\&+\big(-2 a_{n-1,n-1}-(n-2)a_{n-1,n-1}\big)x^{n-1}\\&=
1+\big(2a_{2,n-1}-2a_{1,n-1}\big)x+\sum_{k=3}^{n-1}
x^{k-1}\Big(ka_{k,n-1}-ka_{k-1,n-1}\Big)+\big(-na_{n-1,n-1} \big)x^{n-1}\\&=1+\sum_{k=2}^{n-1}
x^{k-1}\Big(ka_{k,n-1}-ka_{k-1,n-1}\Big)+\big(-n(-1)^{n}(n-1)!\big)x^{n-1}\\&
=1+\sum_{k=2}^{n-1}
x^{k-1}\Big(ka_{k,n-1}-ka_{k-1,n-1}\Big)+(-1)^{n+1}n!x^{n-1} \\&=a_{1,n}+\sum_{k=2}^{n-1}
x^{k-1}\Big(ka_{k,n-1}-ka_{k-1,n-1}\Big)+a_{n,n}x^{n-1}=\sum_{k=1}^{n}a_{k,n}x^{k-1},\end{align*}

since
\begin{align*}
&ka_{k,n-1}-ka_{k-1,n-1}= k \Big\{\sum_{j=0}^k (-1)^{j+1} {k\choose j}j^{n-1} -\sum_{j=0}^{k-1} (-1)^{j+1} {k-1\choose j}j^{n-1}\Big\}\\&=
k\Big\{\sum_{j=0}^{k-1} (-1)^{j+1}j^{n-1}\Big[{k\choose j}-{k-1\choose j}\Big]+ (-1)^{k+1}k^{n-1}\Big\} \\&=
k\Big\{\sum_{j=0}^{k-1} (-1)^{j+1}j^{n-1}{k-1\choose j}\Big[\frac{k}{k-j}-1\Big]+ (-1)^{k+1}k^{n-1}\Big\}\\& =
k\Big\{\sum_{j=0}^{k-1} (-1)^{j+1}j^{n-1}{k-1\choose j}\frac{j}{k-j}+ (-1)^{k+1}k^{n-1}\Big\}\\&=k\Big\{\sum_{j=0}^{k-1} (-1)^{j+1}j^{n}\frac{(k-1)!}{j!(k-j)!}+ (-1)^{k+1}k^{n-1}\Big\}\\&=\sum_{j=0}^{k-1} (-1)^{j+1}j^{n}{k \choose j}+ (-1)^{k+1}k^{n}=\sum_{j=0}^{k} (-1)^{j+1}j^{n}{k \choose j}=a_{k,n}.\end{align*}
This completes the proof.

\hfill $\Box$

\bigskip
\noindent
\begin{corollary}
The following formula holds
$$\frac{c_n(x)}{x}-1=\sum_{k=2}^{n}\frac {a_{k,n}}{k}\,x^{k-1}.$$
\end{corollary}

{\it Proof.} Write
\begin{align*}
&\frac{c_n(x)}{x}= (1-x)\Big\{\sum_{k=1}^{n-1}a_{k,n-1}x^{k-1}\Big\}=
\sum_{k=1}^{n-1}a_{k,n-1}x^{k-1}-
\sum_{k=1}^{n-1}a_{k,n-1}x^{k}\\&=1+\sum_{k=2}^{n-1}a_{k,n-1}x^{k-1}
-
\sum_{k=1}^{n-1}a_{k,n-1}x^{k}=1+\sum_{k=1}^{n-2}
a_{k+1,n-1}x^{k}-
\sum_{k=1}^{n-1}a_{k,n-1}x^{k} \\&=1+\sum_{k=1}^{n-2}\big(a_{k+1,n-1}-a_{k,n-1}\big)x^{k}-a_{n-1,n-1}x^{n-1}=1+\sum_{k=1}^{n-2}\frac {a_{k+1,n}}{k+1}\,x^{k}+(-1)^{n-1}+(n-1)!x^{n-1},\end{align*}
since, from the last calculation above
$$a_{k+1,n-1}-a_{k,n-1}=\frac {a_{k+1,n}}{k+1}.$$
Using now Remark \ref{1} (ii), we get
\begin{align*}
&\frac{c_n(x)}{x}=1+\sum_{k=1}^{n-2}\frac {a_{k+1,n}}{k+1}\,x^{k}+(-1)^{n+1}(n-1)!x^{n-1}\\&=1+\sum_{k=2}^{n-1}\frac {a_{k,n}}{k}\,x^{k-1}+\frac {(-1)^{n+1}n!}{n}x^{n-1}=1+\sum_{k=2}^{n}\frac {a_{k,n}}{k}\,x^{k-1}.
\end{align*}
Thus we have obtained
$$\frac{c_n(x)}{x}-1=\sum_{k=2}^{n}\frac {a_{k,n}}{k}\,x^{k-1},$$
as claimed.

\hfill $\Box$

\end{document}